\documentclass[a4paper,twoside]{article}
\usepackage{amssymb}
\usepackage{amsmath}
\usepackage[utf8]{inputenc}
\usepackage[T1]{fontenc}
\usepackage[final]{lpretex}
\usepackage[all]{xy}
\usepackage{amscd}
\usepackage{title}

\DocVersion[Tue Aug 7 2007]{1}
\DocVersion[Tue Nov 6 2007]{2}

\EndOfDump
\usepackage{htmlref}
\usepackage{url}
\newcommand\Kscheme[2][]{K_0^{#1}(\mathrm{Spc}_{#2})}
\newcommand\Kschemecompl[2][]{\widehat{K}_0^{#1}(\mathrm{Spc}_{#2})}
\newcommand\Kstack[2][]{K_0^{#1}(\mathrm{Stck}_{#2})}
\newcommand\Kcoh[1]{K_0(\mathrm{Coh}_{#1})}
\newcommand\Kcompl[1]{\widehat K_0(\mathrm{Coh}_{#1})}
\newcommand\Gcoh[1]{L_0(\mathrm{Coh}_{#1})}
\newcommand\Gcompl[1]{\widehat L_0(\mathrm{Coh}_{#1})}
\newcommand\Gcplx{L_0(\mathrm{Cplx})}
\newcommand\Gcmpl{\widehat L_0(\mathrm{Cplx})}
\newcommand\Gab[1][]{L_0(\mathrm{Ab}_{#1})}
\newcommand\Kspol[1]{\overline{K}_0^{\mathrm{pol}}(\mathrm{Spc}_{#1})}
\newcommand\Kcpol[1]{\overline{K}_0^{\mathrm{pol}}(\mathrm{Coh}_{#1})}
\newcommand\Lclass{\mathbb{L}}
\newcommand\Bclass{\mathrm{B}}
\newcommand\NS{\mathrm{NS}}

\newcommand\Zar{\Cal{Z}\mathrm{ar}}
\newcommand\uwt{\overline w}
\newcommand\cnt{\symb{cnt}}
\newcommand\Fil{\symb{Fil}}

\newcommand\Gal{\symb{Gal}}
\newcommand\Bun{\frak{Bun}}
\newcommand\hco{{:}}

\begin{document}
\title{The Grothendieck group of algebraic stacks} 

\author{Torsten Ekedahl} 

\address{Matematiska institutionen\\
Stockholms universitet\\ 
SE-106 91 Stockholm\\ 
Sweden}
\email{teke@math.su.se} 

%% 14E08 Rationality questions
%% 14A20 Generalizations (algebraic spaces, stacks)
%% 14G15 Finite ground fields
%% 14F22 Brauer groups of schemes
%% 14J10 Families, moduli, classification: algebraic theory
%% 14F20 Étale and other Grothendieck topologies and cohomologies
\begin{abstract}
We introduce a Grothendieck group of algebraic stacks (with affine stabilisers)
analogous to the Grothendieck group of algebraic varieties. We then identify it
with a certain localisation of the Grothendieck group of algebraic
varieties. Several invariants of elements in this group are discussed. The most
important is an extension of the Euler characteristic (of cohomology with
compact support) but in characteristic zero we introduce invariants which are
able to distinguish between classes with the same Euler characteristic. These
invariants are actually defined on the completed localised Grothendieck ring of
varieties used in motivic integration. In particular we show that there are
$\PSL_n$-torsors of varieties whose class in the completed localised
Grothendieck ring of varieties is not the product of the class of the base and
the class of $\PSL_n$.
\end{abstract}
\maketitle

Counting points of varieties over finite fields is both a well-established
heuristic method to guess more precise (cohomological) statements as well as a
way to obtain actual information on the $\ell$-adic Galois representation Euler
characteristics. Particularly in the latter case what has often been considered
(starting at least with \cite{harder74}) are moduli varieties (or usually more
properly stacks) where one when counting constantly divides by the order of some
automorphism group. There is a more direct approach to the cohomology of the
moduli spaces considered by Harder-Narasimhan given by Atiyah-Bott
(\cite{atiyah83::yang+mills+rieman}). Though completely different in basic
approach it is surprisingly similar in the logic of the recursion which
inductively analyses the strata of a stratification. The main difference in
execution is that point counting allows one to just sum up over the strata
and divide by the order of a group when one is taking quotients by it. That one
can sum over the strata is essentially because one is computing Euler
characteristics of cohomology instead of getting hold of the actual cohomology
(though a clever part of the Harder-Narasimhan argument is that one can compute
the Betti numbers from the Galois representation Euler characteristic in many
cases by using Deligne's fundamental weight result). The second part is somewhat
different as the division by the group order corresponds to a spectral sequence
involving equivariant cohomology. As this spectral sequence will in general have
infinite total dimension it is not possible to directly compute Euler
characteristics from it.

This article started off as an attempt to set things up so that the geometrical
aspects of the point counting of \cite{bergstroem06::point} could be made more
visible. For that one might want to go further than just making a cohomological
argument as the geometric reasons are usually even of a ``motivic''
character.\footnote{I use quotes as I do not use ``motivic'' in any technically
precise way.} Point counting of varieties can often be made motivic by working
in a Grothendieck ring of varieties and the kind of point counting that is
involved here should be made motivic by working in a Grothendieck ring of
algebraic stacks. In such a ring the classes of many algebraic groups will be
invertible so that one may indeed divide by them. I eventually found out that
this idea had already appeared to several people (cf.,
\cite{behrend07,joyce07::motiv+artin,toen05::groth+artin}). I believe however
that we manage to navigate the Scylla and Charibdis of these precedents to the
point that the introduction of a new Grothendieck group of algebraic stacks can
be argued (see the remark after Theorem \ref{Localisation theorem} for an
extended discussion of the relations between the present approach and these
others). A more precise description of the contents is as follows.

The article starts with introducing a Grothendieck ring $\Kstack{\k}$ of
algebraic stacks (and some natural quotients of it) over a field $\k$ and gives
a number of situations when fibrations are multiplicative with respect to the
classes of the stacks in this ring. It is then shown that this Grothendieck ring
is actually a specific localisation of the Grothendieck ring $\Kscheme{\k}$ of
algebraic varieties. This localisation result is then used to give an Euler
characteristic of any algebraic stack (of finite type and with affine
stabilisers) in a suitably completed Grothendieck ring of mixed Galois
representations or mixed Hodge structures. It is further shown that the image
lies in a subring of this completion of elements of ``polynomial growth''. This
subring is given a stronger topology than the completion topology. The main
point about this stronger topology is that over a finite field the trace of the
Frobenius on it is well-defined and continuous. It is then shown that this trace
applied to the class of an algebraic stack is equal to the point count of the
stack. Inspired by \cite{behrend07} these results are extended to certain
algebraic stacks only locally of finite type. The Euler characteristic (for
cohomology with compact support) has been defined by extension from the spatial
case but we give the expected comparison with the recent definition of Laszlo
and Olsson (\cite{laszlo08::artin+ii}) of cohomology with compact support for
algebraic stacks.

Even though the point counting focusses the attention on the appropriate Euler
characteristic of classes in $\Kscheme{\k}$ and its variants there are other
invariants such as those of \cite{larsen04::ration} and \cite{poonen02::groth}
that detect elements not detected by the Euler characteristic. The invariant of
\cite{larsen04::ration} does not extend to even $\Kscheme{\k}[\Lclass^{-1}]$,
where $\Lclass$ is the class of $\A^1$, as it in fact vanishes on $\Lclass$ and
hence cannot be used for our purposes. Over the complex numbers we show how to
extend Poonen's invariant to the completion of $\Kscheme{\k}$ and also define
some other related invariants. We give two applications of these invariants. One
is the generalisation of Poonen's example of zero-divisors in $\Kscheme{\k}$ to
both its completion $\Kschemecompl{\k}$ and to $\Kstack{\k}$. The other is
showing that contrary to the case of $\GL_n$-torsors we do not have
multiplicativity for the classes of $\PSL_n$-torsors in $\Kscheme{\k}$ and its
variants. This is pertinent, in the case $G=\PSL_n$, to the question posed by
Behrend and Dhillon (cf., \cite[Rmk 2.10]{behrend07}) on the relation between
$\Kschemecompl{\k}$ and $\Kschemecompl[G]{\k}$. A more substantial use of these
invariants will however be made in \cite{ekedahl08::stack} where we shall show
the non-triviality of the classes of the classifying stacks of certain finite
groups.

At the end we make a sample computation for the class of the algebraic stack of
smooth proper genus $1$ curves with a polarisation of degree $3$.
%;;%Stacks vs spaces
\begin{section}{Stacks vs spaces}

Recall that $\Kscheme{\k}$ is the Grothendieck group spanned by classes $\{X\}$
of algebraic $\k$-spaces $X$ of finite type with the relations that $\{X\}$ only
depends on the isomorphism class of $X$ and $\{X\}=\{Y\}+\{U\}$ where $Y$ is a
closed subscheme of $X$ and $U$ its complement. (Usually one considers only
schemes but that gives the same group, our choice is motivated by the further
consideration of algebraic stacks which are most closely related to algebraic
spaces.) By analogy we let $\Kstack{\k}$ be the Grothendieck group spanned by
classes $\{\sX\}$ of algebraic stacks $\sX$ of finite type over the field $\k$
all of whose automorphism group schemes are affine (we shall assume that the
automorphism group schemes of all our stacks of this type without further
mention) with relations
\begin{itemize}
\item $\{\sX\}$ depends only on the isomorphism class of $\sX$,

\item $\{\sX\}=\{\sY\}+\{\sU\}$, where $\sY$ is a closed substack of $\sX$ and
$\sU$ its complement and

\item $\{\sE\}=\{\sX\times\A^n\}$, where $\sE \to \sX$ is a vector bundle of
constant rank $n$.\footnote{In the case of $\Kscheme{\k}$ this relation follows
from the first two.}
\end{itemize}
As in the case of $\Kscheme{\k}$ we have a ring structure on $\Kstack{\k}$ given
by $\{\sX\}\{\sY\}=\{\sX\times\sY\}$. Putting $\Lclass:=\{\A^1\}$, the class of
the affine line, the last condition above can then be rewritten as
$\{\sE\}=\Lclass^n\{\sX\}$. If $\sG$ is a class of connected $\k$-group schemes
of finite type, then we define $\Kstack[\sG]{\k}$ as the quotient of
$\Kstack{\k}$ by the relations $\{Y\}=\{G\}\{X\}$ for every $G \in \sG$ and
every $G$-torsor $X \to Y$, where $X$ and $Y$ are \emph{spaces}. We have a
unique ring structure on $\Kstack[\sG]{\k}$ making the quotient map
$\Kstack{\k}\to\Kstack[\sG]{\k}$ a ring homomorphism. Note also that if $\sG
\subseteq \sG'$, then the quotient map $\Kstack{\k}\to\Kstack[\sG']{\k}$ factors
through $\Kstack{\k}\to\Kstack[\sG]{\k}$. By considering an algebraic space as
an algebraic stack we of course also get a natural map $\Kscheme{\k} \to
\Kstack{\k}$ as well as maps $\Kscheme[\sG]{\k} \to \Kstack[\sG]{\k}$ with the
obvious definition of $\Kscheme[\sG]{\k}$. We start by collecting some basic
properties of $\Kstack[\sG]{\k}$ in the following proposition.
%;;%%Proposition Base formulas (iii)Classifying inverse (iva) Isotorsor (vi)
%% Subgroup formula (iv) Fibration
\begin{proposition}\label{Base formulas}
\part[i] We have
$\{\GL_n\}=(\Lclass^n-1)(\Lclass^n-\Lclass)\cdots(\Lclass^n-\Lclass^{n-1})$.

\part[ii] If $\sX \to \sY$ is a $\GL_n$-torsor of algebraic $\k$-stacks of finite
type, then $\{\sX\}=\{\GL_n\}\{\sY\}\in \Kstack{\k}$.

\part[iiia] If $G \in \sG$ and $\sX \to \sY$ is a $G$-torsor of algebraic
$\k$-stacks of finite type, then $\{\sX\}=\{G\}\{\sY\} \in \Kstack[\sG]{\k}$.

\part[iiib] If, $H,N \in \sG$, $1 \to N \to G \to H \to 1$ is an extension of
algebraic groups and $\sX \to \sY$ is a $G$-torsor, then $\{\sX\}=\{G\}\{\sY\}$
in $\Kstack[\sG]{\k}$.

\part[iii]\label{Classifying inverse} If $G \in \sG$ then $\{G\}\{\Bclass G\}=1
\in \Kstack[\sG]{\k}$.

\part[iv]\label{Fibration} If $G \in \sG$, $F$ is a $G$-space, $\sX \to \sY$ is
a $G$-torsor of algebraic $\k$-stacks of finite type and $\sZ \to \sY$ the
$F$-fibration associated to $\sX \to \sY$ and the $G$-action on $F$, then
$\{\sZ\}=\{F\}\{\sY\}$.

\part[iva]\label{Isotorsor} If $G \in \sG$, $\sX \to \sY$ and $\sX' \to \sY$ are
$G$-torsors and $\sZ \to \sY$ the stack of isomorphisms $\sX \riso \sX'$, then
$\{\sZ\}=\{G\}\{\sY\}$ in $\Kstack[\sG]{\k}$.

\part[v] If $G \in \sG$, $F$ is a $G$-space, $H$ is an algebraic group, $H \to G$ a
morphism of algebraic groups, $\sX \to \sY$ an $H$-torsor of
algebraic $\k$-stacks of finite type and $\sZ \to \sY$ the $F$-fibration
associated to $\sX \to \sY$ and the $H$-action on $F$ (given by its $G$-action
and the homomorphism $H \to G$), then
$\{\sZ\}=\{F\}\{\sY\}$.

\part[vi]\label{Subgroup formula}If $G \in \sG$ and $H \hookrightarrow G$ is a
subgroup scheme, then $\{\Bclass H\}=\{G/H\}\{\Bclass G\} \in \Kstack[\sG]{\k}$.
\begin{proof}
To start with \DHrefpart{i} we note that $\GL_n$ is isomorphic to the space of
bases of $\A^n$. The space of $k+1$ linearly independent vectors fibres over the
space of $k$ linearly independent vectors with total space $\A^n$ times the base
minus a vector bundle of rank $k$. The fact that $\{E\}=\Lclass^k\{X\}$ for a
vector bundle $E \to X$ of rank $k$ implies \DHrefpart{i} by induction.

As for \DHrefpart{ii} its proof is the relative version of the proof of
\DHrefpart{i}; we get a succession of fibrations $\sX=\sX_{n} \to
\sX_{n-1}\to\dots\to \sY$, where $\sX_k$ is the stack of $k$ linearly
independent vectors in the vector bundle associated to $\sX$. Each
$\sX_{k+1}\to\sX_k$ is then the complement of a vector subbundle in a vector
bundle and we use repeatedly the third set of defining relations of
$\Kstack{\k}$ to get that
$\{\sX\}=(\Lclass^n-1)(\Lclass^n-\Lclass)\cdots(\Lclass^n-\Lclass^{n-1})\{\sY\}$.
We then use \DHrefpart{i} to conclude.

To prove \DHrefpart{iiia}, we start by noticing that $\{\GL_n\}$ is invertible
in $\Kstack{\k}$ for all $n$. Indeed, $\Spec\k \to \Bclass\GL_n$ is a
$\GL_n$-torsor and hence we have $1=\{\GL_n\}\{\Bclass\GL_n\}$. In particular
for a $\GL_n$-torsor $\sZ \to \sW$ we may write the relation of \DHrefpart{ii}
as $\{\sW\}=\{\Bclass\GL_n\}\{\sZ\}$. 

Continuing, the formula to be proved is additive in decompositions of the base
$\sY$ as the union of a closed substack and its complement. Hence we may assume
that $\sY$ is a global quotient, i.e., there is $\GL_n$-torsor $Y \to \sY$ such
that $Y$ is an algebraic space. Hence we have a $2$-cartesian diagram
\begin{displaymath}
\begin{CD}
X @>>> Y \\
@VVV @VVV\\
\sX @>>> \sY,
\end{CD}
\end{displaymath}
where $X \to Y$ is a $G$-torsor of algebraic spaces and $X \to \sX$ and $Y \to
\sY$ are $\GL_n$-torsors. This gives
$\{\sX\}=\{\Bclass\GL_n\}\{X\}=\{\Bclass\GL_n\}\{G\}\{Y\}=\{G\}\{\sY\}$ which
proves \DHrefpart{iiia}. Then, applying \DHrefpart{iiia} to the $N$-torsor $\sX
\to \sX/N$ and to the $H$-torsor $\sX/N \to \sY$ we get
$\{\sX\}=\{N\}\{H\}\{\sY\}$ but $G \to H$ is an $N$-torsor so that
$\{G\}=\{N\}\{H\}$ which gives \DHrefpart{iiib}. Further, \DHrefpart{iii} then
follows by applying \DHrefpart{iiia} to the $G$-torsor $\Spec\k \to \Bclass G$.

Now, to prove \DHrefpart{iv}, the $F$-fibration associated to $\sX \to \sY$ is
the quotient $F\times_G\sX$. As the quotient map $F\times \sX \to F\times_G\sX$
is a $G$-torsor we get $\{F\times_G\sX\}=\{\Bclass G\}\{F\}\{\sX\}=\{\Bclass
G\}\{F\}\{G\}\{\sY\}=\{F\}\{\sY\}$.

In \DHrefpart{iva}, $\sZ$ is the fibration associated to the $G\times G$-torsor
$\sX\times_{\sY}\sX' \to \sY$ and the $G\times G$-space $G\times G/G$, where $G
\subseteq G\times G$ is the diagonal embedding. Hence, \DHrefpart{iva} follows
from \DHrefpart{iv} and \DHrefpart{iiib}.

Turning to \DHrefpart{v} we note that the associated $F$-fibration,
$F\times_H\sX$, is isomorphic to the $F$-fibration associated to the $G$-torsor
$G\times_H\sX \to \sY$. Hence, \DHrefpart{v} follows from \DHrefpart{iv} applied
to this $G$-torsor.

Finally, \DHrefpart{vi} follows directly from \DHrefpart{iv} applied to the
$H$-torsor $\Spec\k \to \Bclass G$ and the $G$-space $G/H$, using that $\Bclass
H= G/H\times_G\Bclass G$.
\end{proof}
\end{proposition}
Using the results just obtained we can obtain a precise relation between
$\Kscheme{\k}$ and $\Kstack{\k}$. Indeed, we let $\Kscheme[\sG]{\k}'$ be the ring
obtained from $\Kscheme[\sG]{\k}$ by inverting $\Lclass$ as well as $(\Lclass^n-1)$
for all $n>0$.
%;;%%Theorem Localisation theorem
\begin{theorem}\label{Localisation theorem}
For any $\sG$, the natural map $\Kscheme[\sG]{\k}\to \Kstack[\sG]{\k}$ induces
an isomorphism $\Kscheme[\sG]{\k}'\to \Kstack[\sG]{\k}$.
\begin{proof}
By Proposition \ref{Base formulas} we have that
$\{\GL_n\}=(\Lclass^n-1)(\Lclass^n-\Lclass)\cdots(\Lclass^n-\Lclass^{n-1})$ is
invertible in $\Kstack[\sG]{\k}$ which implies that we do indeed get a
factorisation $\Kscheme[\sG]{\k}'\to \Kstack[\sG]{\k}$ and we shall proceed by
defining a map in the other direction. If an algebraic stack $\sX$ (of finite
type and with affine stabilisers) is a global quotient $\sX=[X/\GL_n]$, then we
would like to define $\{\sX\}$ in $\Kscheme[\sG]{\k}'$ as $\{X\}/\{\GL_n\}$. Our
localisation ensures that $\{\GL_n\}$ is indeed invertible. If we present
$\sX$ in a second way as a global quotient $\sX=[X'/\GL_{n'}]$ then we may
construct a $2$-cartesian diagram
\begin{displaymath}
\begin{CD}
X'' @>>> X\\
@VVV @VVV\\
X' @>>> \sX
\end{CD}
\end{displaymath}
and as $X'' \to X$ is a $\GL_{n'}$-torsor and $X'' \to X'$ a $\GL_n$-torsor we
get $\{\GL_{n'}\}\{X\}=\{X''\}=\{\GL_n\}\{X'\}$ which gives the independence. In
the general case we stratify a stack $\sX$ by global quotients and sum up. To
show well-definedness we may assume that $\sX$ is a global quotient and we get
what we want by the independence just proved.
\end{proof}
\end{theorem}
\begin{remark}
As mentioned in the introduction this theorem is already implicitly, though it
seems not explicitly, in the literature. In \cite{behrend07} the authors
associate to certain algebraic stacks (including those of finite type with
affine stabilisers) an element in $\Kcompl{\k}$. An analysis of the proof shows
that if one restricts oneself to finite type stacks then one need only invert
$\Lclass$ and $\Lclass^n-1$ in $\Kscheme{\k}$ for the construction to work. In
any case the crucial point, that $\{X\}\{\GL_n\}^{-1}$ only depends on the stack
$[X/\GL_n]$ and not on the actual action of $\GL_n$ on the algebraic space $X$
is to be found in \cite{behrend07}. This is also \cite[Prop.\
4.8]{joyce07::motiv+artin} (though formally only in
$\Kstack{\k}\Tensor\Q$). Finally, one can obtain the map $\Kstack{\k} \to
\Kscheme{\k}'$ by including algebraic stacks in the larger category of
$n$-stacks and using Toën's $n$-analogue, \cite[Thm 3.10]{toen05::groth+artin},
of the theorem. However, it seems worthwhile to give an explicit definition of a
Grothendieck group of algebraic stacks whose relations are simpler than a
restriction of the relations of Toën to the case of $1$-stacks. Also giving a
direct proof of the theorem is arguably more convenient for the reader than a
piecing together of the arguments found in \cite{behrend07}and
\cite{joyce07::motiv+artin}. Note also that parts of Proposition \ref{Base
formulas} are also implicitly in the quoted sources (or analogues in the case of
\cite{behrend07}).
\end{remark}
It is not at all clear how drastic the passage from $\Kscheme{\k}$ to
$\Kscheme{\k}'$ is, we don't for instance know anything about the kernel of the
localisation map $\Kscheme{\k} \to \Kscheme{\k}'$. The following proposition
gives a small amount of information.
%;;%%%Proposition
\begin{proposition}
Suppose that $\varphi(x)$ and $\rho(x)$ are integer polynomials. If $\varphi(\Lclass)$ divides
$\rho(\Lclass)$ in $\Kscheme{\k}$, then $\varphi(x)$ divides $\rho(x)$. In particular the
only integer polynomials in $\Lclass$ that are invertible in $\Kscheme{\k}$ resp.\
$\Kstack{\k}$ are $\pm1$ resp.\ products of\/ $\Lclass$ and cyclotomic polynomials.
\begin{proof}
Assume that $\rho(\Lclass)=\varphi(\Lclass)x$ for some $x \in
\Kscheme{\k}$. Applying $\chi_c$ (whose definition is recalled below) we get an
equation $\rho(q)=\varphi(q)y$ in $\Kcoh{\k}$. We now argue by induction on the
degree of $\rho$. Let $y'$ be the term of highest weight part of $y$ and
$a_nq^n$ the highest degree term of $\varphi(q)$. Then $a_nq^ny'$ is the highest
weight part of $\rho(q)$ and is hence of the form $b_mq^m$, where $\rho$ is of
degree $m$. As $q$ is invertible we get $a_nq^{n-m}y'=b_m$. As $\Kcoh{\k}$ is a
free abelian group (free on the irreducible Galois representations or rational
Hodge structures) we get that $a_n$ divides $b_m$. Thus we get
$\rho-(a_n/b_m)q^{n-m}\varphi(q)=\varphi(q)(y-(a_n/b_m)q^{n-m})$ and the
induction assumption implies that $\varphi$ divides $\rho$. The rest of the
proposition follows from the first part (using that
$\Kstack{\k}=\Kscheme{\k}'$).
\end{proof}
\end{proposition}
Now, let $\Zar=\Zar_{\k}$ be the class of connected finite type group schemes
for which torsors over any (finitely generated) extension field of $\k$ are
trivial. We have of course that $\GL_n, \SL_n \in \Zar$ but the following
proposition provides more examples.
%;;%%Proposition Zariski locally trivial
\begin{proposition}\label{Zariski locally trivial}
\part[i] We have that $\Kstack{\k}=\Kstack[\Zar]{\k}$.

\part[ii] Let $L$ be a finite dimensional $\k$-algebra. Then $L^*$, the group
scheme of invertible elements of $L$, belongs to $\Zar$.

\part $\Zar$ is closed under extensions.
\begin{proof}
For \DHrefpart{i} using Theorem \ref{Localisation theorem} it is enough to show
that for any $G \in \Zar$ and any $G$-torsor $X \to Y$ of algebraic $\k$-spaces
of finite type we have that $\{X\}=\{G\}\{Y\} \in \Kscheme{\k}$. It is enough to
find a stratification of $Y$ such that the torsor is trivial over any
stratum. We do this by Noetherian induction using that over each generic point of
$Y$ the torsor is trivial by assumption and the trivialisation can be extended
to an open subset of $Y$.

For a $\k$-field $K$, $L^*$-torsors over $K$ correspond to $L_K$-modules locally
free of rank $1$ and such modules are indeed free which proves \DHrefpart{ii}.

The fact that $\Zar$ is closed under extensions is clear.
\end{proof}
\end{proposition}
\begin{remark}
Over an algebraically closed field the groups in $\Zar$ are precisely the
special groups in the sense of \scite[Exp.\ 1]{58::semin+c}. Indeed, a slight
modification of the proof of \lcite[Thm 1]{58::semin+c} shows that a group
$G\in\Zar$ is affine (and connected) and to show that it is special it suffices
by \lcite[Thm 2]{58::semin+c} to show that $\GL_n \to \GL_n/G$ is Zariski
locally trivial but as observed after that theorem this follows from the
existence of a rational section which in turn it has as $G \in \Zar$. The third
part of the proposition is, again in the algebraically closed case, \lcite[Lemme
6]{58::semin+c} while the second is in \cite{joyce07::motiv+artin}. Note that
while the proof of these results in the non-algebraically closed case are the
same as for the algebraically closed case some care has to be taken. For
instance \lcite[Prop.\ 14]{58::semin+c}, that solvable connected groups are
special, is not true in general. There are in fact tori over suitable fields
that are not special (or in $\Zar$) and over non-perfect fields there are also
non-special connected unipotent algebraic groups.
\end{remark}
\end{section}
%;;%Section Euler characteristics and point counting
\begin{section}{Euler characteristics and point counting}

We start by recalling the definition of the Euler characteristic of elements of
$\Kscheme{\k}$ and $\Kschemecompl{\k}$. Recall that the Euler characteristic of
cohomology with compact support is additive over a decomposition of an algebraic
space into a closed subspace and its complement. This means that this Euler
characteristic gives a map from $\Kscheme{\k}$. We want the recipient group to
be an appropriate Grothendieck group of cohomology groups with extra
structure. When $\k=\C$ one such choice is the Grothendieck group of
(polarisable) mixed Hodge structures (and when $\k=\R$ one could consider mixed
Hodge structures with a complex conjugation). For a general field we would like
to use actions of Galois groups of finitely generated subfields of $\k$. To be
able to get the needed functoriality as well as to have a notion of mixedness we
need to spread out even further. Hence we shall be considering the category that
is the direct limit of the category of local $\Q_{\ell}$-systems on (spectra of)
finitely generated subrings of $\k$. In concrete Galois-theoretic terms this
means the following: Fixing an algebraic closure $\bk$ of $\k$ we can consider,
for each finitely generated subring $S$ of $\k$, the fundamental group
$\pi_1(\Spec S,\overline s)$, where $\overline s$ is the composite $\Spec\bk \to
\Spec\k\to \Spec S$. An object in the category is then a continuous
finite-dimensional $\Q_{\ell}$-representation of $\pi_1(\Spec S,\overline s)$
for some $S$. A morphism between two such representations associated to $R$ and
$S$ is then a choice of a $T$ containing both $R$ and $S$ and a map between the
induced $\pi_1(\Spec T,t)$-representations. Such a representation is
\emph{mixed} if, after possibly extending $S$, it is a successive extension of
pure representations. A representation is \emph{pure} (of weight $n$) if for
every closed point $x$ of $\Spec S$, the eigenvalues of the (geometric)
Frobenius map associated to $x$ are algebraic numbers all of whose (archimedean)
absolute values are equal to $q^{n/2}$, where $q=|\k(x)|$. Note that when $S$ is
normal (which can be insured by extending $S$) the map $\symb{Gal}(\overline
K/K) \to \pi_1(\Spec S,s)$, where $K$ is the fraction field of $S$ and
$\overline K$ the algebraic closure of $K$ in $\bk$, is surjective. Hence
$\pi_1(\Spec S,s)$-representations can be thought of as $\symb{Gal}(\overline
K/K)$-representations for which the kernel of $\symb{Gal}(\overline K/K) \to
\pi_1(\Spec S,s)$ acts trivially. In this way our objects can indeed be thought
of as Galois representations for finitely generated subfields of $\k$. When $\k$
is algebraically closed things can be made even more concrete: One chooses a
finitely generated subring $S$ of $\k$ and considers representations of open
subgroups of $\pi_1(\Spec S,\overline s)$ where morphisms are only required to
commute with the elements of some open subgroup. In the general case one does
not allow all open subgroups but only some. If one does not wish to take this
direct limit of categories of mixed representations over finitely generated
subfields one may instead define $\Kcoh{\k}$ as follows:
\begin{itemize}
\item If $\k$ is a finitely generated field, then $\Kcoh{\k}$ is the
Grothendieck group of the category of mixed
$\Gal(\bk/\k)$-representations, i.e., $\Q_\ell$-representations $\Gal(\bk/\k)
\to \GL(V)$ which factor through a quotient $\Gal(\bk/\k) \to \pi_1(S,\bk)$,
where $S$ is a normal integral scheme with $\k$ as function field of finite type
over $\Z$ and which are mixed. This group does not depend on $\bk$ and is
functorial for field inclusions.

\item For a general field $\Kcoh{\k}:=\dli_K\Kcoh{K}$, where $K$ runs over
finitely generated subfields of $\k$.
\end{itemize}

In any case we shall use $\Kcoh{\k}$ to denote the Grothendieck group of either
mixed Galois $\Q_{\ell}$-representations as above or when $\k=\C$ ($\k=\R$) the
mixed polarisable Hodge structures (resp.\ those provided with a complex
conjugation). The additivity of the Euler characteristic (of cohomology with
compact support) can then be formulated as saying that it gives a ring
homomorphism \map{\chi_c}{\Kscheme{\k}}{\Kcoh{\k}} (where the ring structure on
$\Kcoh{\k}$ is induced by the tensor product of representations). Note also that
the Euler characteristic is multiplicative for torsors over a connected group
scheme (indeed all higher direct images of a constant sheaf are constant as they
become constant after a base change with connected fibres) so that $\chi_c$
factors over a \map{\chi_c}{\Kscheme[\sG]{\k}}{\Kcoh{\k}} for an arbitrary
$\sG$.

Now, $\chi_c(\Lclass)$ is invertible but $\chi_c(\Lclass^n-1)$ is not. In order
for us to be able to use Theorem \ref{Localisation theorem} to extend $\chi_c$
to $\Kstack{\k}$ we could of course just invert the
$\chi_c(\Lclass^n-1)$. However, we shall follow the established road of making a
completion instead. For this note that as every cohomology object has finite
length and the simple objects are pure of some weight, $\Kcoh{\k}$ is graded by
weight. Hence, the completion, which completes with respect to the weight
filtration in the direction of negative weights is very benign in that it just
replaces finite sums with sums unbounded in the negative direction. We shall use
$\Kcompl{\k}$ to denote this completion. Putting $q:=\chi_c(\Lclass)$ we then
have $\chi_c(\Lclass^n-1)=q^n-1$ which is invertible with inverse
$q^{-n}+q^{-2n}+\cdots$ so that we do indeed get an extension
\map{\chi_c}{\Kstack{\k}}{\Kcompl{\k}}.

We now want to compare this with the process of point counting. Hence, if $\F$
is a finite field we have a map \map{|-|}{\Kstack{\F}}{\Q} by associating to an
algebraic stack (of finite type) $\sX$ the mass of the essentially finite
groupoid $\sX(\F)$ (the \Definition{mass} of an essentially finite groupoid is
the sum over its isomorphism classes of $1$ over the order of the automorphism
groups). On the other hand, for an element $\sum_{i\le N}a_i$ of $\Kcompl{\F}$
(with $a_i$ pure of weight $i$) one can consider the sum $\sum_{i\le N}\Tr(a_i)$
where \map{\Tr}{\Kcoh{\F}}{\bQ} is the map that takes a Galois representation to
the trace of the action of the Frobenius element (where we have fixed once and
for all an embedding of the algebraic numbers of $\Q_\ell$ into $\bQ$). In
general this sum will not be convergent and one needs to add some conditions on
a sum in order to be allowed to take its trace. In \cite{behrend07} the authors
use what seems essentially to be the weakest possible condition namely that the
sum is absolutely convergent. However, this seems to me to be too \emph{ad hoc},
for instance possibly depending on the precise finite field that one is
considering and also makes sense only when the base field is finite. I instead
propose the following approach:

Given a $k$-scheme $X$ of finite type and an integer $n$ we define $w_n(X)$ to
be $\sum_iw_n(H^i_c(X))$, where $w_n(H^i_c(X))$ is the dimension of the part of
$H^i_c(X)$ of weight $n$ (note that in the case when $\k$ equals $\C$ or $\R$
this is not ambiguous as the weight filtration on the classical cohomology is
compatible with that on \'etale cohomology).  We then define, for $x \in
\Kscheme{\k}$ $w'_n(x)$ to be the minimum of $\sum_kw_n(X_k)$ where
$\sum_k\pm\{X_k\}$ runs over all representatives of $x$ as a signed sum of
classes of schemes. Note that $w_{n+2}(X\times \A^1)=w_{n}(X)$ and thus
$w'_{n+2}(x\Lclass) \le w'_{n}(x)$. This means that we may put
$\uwt_n(x):=\lim_{k\to\infty}w'_{n+2k}(x\Lclass^k)$.  We now extend $\uwt_n$ to
$\Kscheme{\k}[\Lclass^{-1}]$ by putting $\uwt_n(x\Lclass^{-k}):=\uwt_{n+2k}(x)$
(this is consistent as by construction $\uwt_{n+2}(x\Lclass)=\uwt_{n}(x)$ for $x
\in \Kscheme{\k}$). We clearly have subadditivity for $\uwt_n$, $\uwt_n(x\pm
y)\le \uwt_n(x)+\uwt_n(y)$, as well as submultiplicativity, $\uwt_n(xy)\le
\sum_{i+j=n}\uwt_i(x)\uwt_j(y)$, by the K\"unneth formula. The next step is to
complete $\Kscheme{\k}[\Lclass^{-1}]$ by adding equivalence classes of Cauchy
sequences $x=(x_i)$ fulfilling two conditions. The first uses the dimension
filtration $\{\Fil^{\le m}\Kscheme{\k}[\Lclass^{-1}]\}$ where $\Fil^{\le
m}\Kscheme{\k}[\Lclass^{-1}]$ is spanned by classes $\{X\}\Lclass^{-k}$ with
$\dim(X)-k\le m$ and $\dim(x)$ is the least (possibly equal to $-\infty$) number
such that $x \in \Fil^{\le m}\Kscheme{\k}[\Lclass^{-1}]$. We then, to begin
with, demand that $\dim(x_i-x_j) \to -\infty$ when $i,j \to \infty$. We further
impose that there be constants $C$ and $D$ both independent of $i$ such that
$\uwt_n(x_i)\le C|n|^d+D$ for all $n$. We shall say that such a sequence is
\Definition{convergent of uniform polynomial growth}.  Note that by results of
Deligne (cf., \cite{deligne80::la+weil,deligne74::theor+hodge+iii}) we have that
$\uwt_n(x) = 0$ if $2\dim(x)<n$ and thus (by subadditivity) for any $n$,
$\uwt_n(x_i)$ is eventually constant and we may unambigously define $\uwt_n(x)$
as that constant value. As $(x_i)$ is of uniformly polynomial growth we also get
that $\uwt_n(x)$ is of polynomial growth.

The condition of polynomial growth can be conveniently formulated in the
following way. For $x \in \Kscheme{\k}[\Lclass^{-1}]$ we put
$\uwt(x,t)=\sum_{n}\uwt_n(x)t^{n}$ considered as an element of $\Z((t^{-1}))$,
Laurent series in $t^{-1}$ with integer coefficients. That a sequence $(x_i)$ is
of uniform means that $\uwt(x_i,t) \le \varphi(t)/(t-1)^d$ for some Laurent
polynomial $\varphi$ and an integer $d$ both independent of $i$ and where
comparison is made coefficient-wise. Furthermore, subadditivity can be
formulated as $\uwt(x\pm y,t)\le \uwt(x,t)+\uwt(x,t)$ and submultiplicativity as
$\uwt(xy,t)\le\uwt(x,t)\uwt(y,t)$. These facts make it obvious that convergent
sequences uniformly of polynomial growth are closed under term-wise addition,
subtraction and multiplication. As those that converge to zero (wrt the
dimension filtration) form an ideal we may define $\Kspol{\k}$ as the quotient
of the convergent sequences uniformly of polynomial growth by those converging
to zero. We can extend the $\uwt_n$ (or equivalently $\uwt(-,t)$) to
$\Kspol{\k}$ and it is clear that convergent sequences uniformly of polynomial
growth of elements of $\Kspol{\k}$ converge in $\Kspol{\k}$. We of course have a
completion map \map{\overline{(-)}}{\Kscheme{\k}[\Lclass^{-1}]}{\Kspol{\k}}. We
also get a continuous map from $\Kspol{\k}$ to the completion
$\Kschemecompl{\k}$ of $\Kscheme{\k}[\Lclass^{-1}]$ in the topology given by the
dimension filtration. By definition a uniformly convergent sequence converges to
zero precisely when it converges to zero in the dimension filtration
topology. Hence the map $\Kspol{\k}\to \Kschemecompl{\k}$ is injective and we
may think of $\Kspol{\k}$ as a subring of $\Kschemecompl{\k}$. Its topology is
however stronger than the subring topology.

In a similar manner we also define a subring $\Kcpol{\k}$ of $\Kcompl{\k}$
consisting of those elements $x$ for which $\uwt_n(x)$ has at most polynomial
growth. Here we define $\uwt_n(x)$ to be the sum $\sum_k|m_k|\dim V_k$, where
$\sum_km_k[V_k]$ is the weight $n$-component of $x$ written as a sum of
non-isomorphic irreducible Galois representations (resp.\ Hodge structures). It
is clear that if $x \in \Kscheme{\k}$, then $\uwt_n(\chi_c(x))\le \uwt_n(x)$ so
that $\chi_c$ extends to a map $\Kspol{\k}\to\Kcpol{\k}$. We want however to
give $\Kcpol{\k}$ a stronger topology than that provided by $\Kcompl{\k}$. Thus
we declare a sequence $(x_i)$ in $\Kcpol{\k}$ be \Definition{convergent of
uniform polynomial growth} if it is convergent in the filtration topology (i.e.,
for every $n$ the weight $n$-components of the $x_i$ is eventually constant) and
$\uwt_n(x_i)$ is uniformly of polynomial growth as above. Then $(x_i)$ converges
in $\Kcompl{\k}$ to an element of $\Kcpol{\k}$. (The topology is defined by
having a subset being closed if it is closed under uniformly convergent
sequences.) Given an element $\sum_nx_n$ in $\Kcpol{\k}$, the partial sums
$\sum_{n\ge -N}x_n$ converge uniformly to the element so that $\Kcoh{\k}$ is
dense in $\Kcpol{\k}$. Furthermore, as $\uwt_n(\chi_c(x)) \le \uwt_n(x)$ we get
a continuous extension \map{\chi_c}{\Kspol{\k}}{\Kcpol{\k}} which is a ring
homomorphism.

Before our next result we need to recall that for an algebraic stack $\sX$
locally of finite type over a finite field $\k$ we can define its
\Definition{point count}, $\cnt(\sX) \in [0,\infty]$, to be the mass of the
groupoid $\sX(\k)$, i.e., the sum $\sum_x1/|\Aut(x)|$, where the sum runs over
the isomorphism classes of $\sX(\k)$ (and $\Aut(x)$ is the group of
automorphisms \emph{defined} over $\k$). Being a sum over non-negative numbers
it is always convergent though it may be $\infty$; when $\sX$ is of finite type
it is a finite sum however. The point count is additive over a decomposition
into a closed substack and its complement and hence gives a ring homomorphism
\map{\cnt}{\Kstack{\k}}{\Q}. Note that its composite with the natural map
$\Kscheme{\k} \to \Kstack{\k}$ is the usual counting map. Finally, we have a map
\map{\Tr}{\Kcoh{\k}}{\bQ} given by associating to any Galois representation the
trace of the Frobenius map and where we once and for all have chosen an
embedding of $\bQ$, the field of complex algebraic numbers, in an algebraic
closure of $\Q_\ell$. We now have the following (expected) compatibility result.
For the second part we define, in analogy with the spatial case, the dimension filtration
$\Fil^{\le m}\Kstack{\k} \subseteq \Kstack{\k}$ to be the subgroup spanned by
the algebraic stacks of dimension $\le m$. We also extend $\Fil^{\le
m}\Kscheme{\k}$ to $\Kspol{\k}$ by taking its closure.
%;;%%Proposition Finite type point counting (ia)completion continuity
\begin{proposition}\label{Finite type point counting}
\part[i] The elements $\Lclass$ and the $\Lclass^n-1$, all $n>0$, of
$\Kspol{\k}$ are invertible so that the completion map extends to a ring
homomorphism $\Kscheme{\k}' \to \Kspol{\k}$ and thus to a map $\Kstack{\k} \to
\Kspol{\k}$ which shall also be denoted $\overline{(-)}$.

\part[ia]\label{completion continuity} The completion map
\map{\overline{(-)}}{\Kstack{\k}}{\Kspol{\k}} takes $\Fil^{\le m}\Kstack{\k}$ to
$\Fil^{\le m}\Kspol{\k}$. In particular it is continuous for the filtration
topologies.

Assume now that $\k$ is a finite field.

\part[ii] For every $x=\sum_{n\le N}x_n \in \Kcpol{\k}$, where $x_n$ is the
component of weight $n$ the sum $\sum_n\Tr(x_n)$ is convergent in $\C$ and gives
a \emph{continuous} extension of $\Tr$ to a ring homomorphism
\map{\Tr}{\Kcpol{\k}}{\C}.

\part[iii] For every $x \in \Kstack{\k}$ we have that
$\cnt(x)=\Tr(\chi_c(\overline{x}))$.

\part[iv] \map{\cnt}{\Kscheme{\k}[\Lclass^{-1}]}{\Q} extends to a continuous
ring homomorphism \map{\cnt}{\Kspol{\k}}{\R}.
\begin{proof}
We have that $\Lclass$ is invertible already in
$\Kscheme{\k}[\Lclass^{-1}]$. Furthermore,
$\Lclass^n-1=\Lclass^n(1-\Lclass^{-n})$ and the series $\sum_{k\ge
-1}\Lclass^{-kn}$ is a convergent sum, uniformly of polynomial growth. This
proves \DHrefpart{i}.

To prove \DHrefpart{ia} we may consider only the case $x=\{\sX\}$ where
$\dim(\sX)\le m$ and by cutting up $\sX$ we may further assume that
$\sX=[X/\GL_n]$ so that $\dim X \le m+n^2$. We then have
$\overline{\{\sX\}}=\{X\}\{\GL_n\}^{-1}$ but $\{\GL_n\}^{-1}$ is a power series
in $\Lclass^{-1}$ starting with $\Lclass^{-n^2}$ so that $\{X\}\{\GL_n\}^{-1}
\in \Fil^{\dim X-n^2}\Kspol{\k}=\Fil^{m}\Kspol{\k}$.

As for \DHrefpart{ii}, if $\k$ has $q$ elements then by definition there is a
$d$ such that $|\Tr(x_n)|\le |n|^dq^{n/2}$ for all $n$ which gives
convergence as $q>1$. Furthermore, if $(x_i)$ is a convergent sequence of
uniform polynomial growth converging in $\Kcpol{\k}$ to $x$, then there are $C$,
$D$ and $d$ such that $\uwt_n{x_i} \le C|n|^d+D$ for all $n$ and $i$ and then 
\begin{displaymath}
|\Tr(x)-\Tr(x_i)| \le \sum_{k\le 2\dim(x-x_i)}2(C|k|^d+D)q^{k/2},
\end{displaymath}
which shows that as $\dim(x-x_i) \to -\infty$ as $i \to \infty$ we have
$\Tr(x_i) \to \Tr(x)$.

Further, as both sides are ring homomorphisms, to verify the equality of
\DHrefpart{iii}, it is enough to check it on $\Kscheme{\k}$ where it is simply
the Lefschetz fixed point formula.

Finally, to prove \DHrefpart{iv} we note that $\Tr\circ \chi_c$ is a continuous
extension of $\cnt$ to $\Kspol{\k}$.
\end{proof}
\end{proposition}
\begin{remark}
\part For the truth of the proposition it would be enough to have $\uwt_n(x)$
grow slower than $\alpha^{-n}$ for all $\alpha > 1$ (or even for some $\alpha <
2$). However, the requirement of polynomial growth seems to me to be the more
natural condition (see however comment below).

\part Note that $\uwt_n(\chi_c(x))$ (for $x \in \Kscheme{\k}$ say) is in general
smaller than $\uwt_n(x)$ because of possible cancellations of irreducible Galois
representations which do not correspond to cancellations of varieties in
$\Kscheme{\k}$. Working with a Grothendieck group of motives instead of
varieties would no doubt get a value of $\uwt_n$ closer to that of
$\Kcpol{\k}$. On the other hand staying in $\Kscheme{\k}$ gives more
information. I also find it quite difficult to imagine that there are natural
examples of elements $x$ in the dimension completion of $\Kscheme{\k}$ for which
$\uwt_n(x)$ does not have polynomial growth while $\uwt_n(\chi_c(x))$ does.

\part The extension of $\chi_c$ to $\Kstack{\k}$ and the relation with point
counting for $\k$ a finite field is done in \cite{toen05::groth+artin} for the
general case of $n$-stacks using the same idea.
\end{remark}
%;;%%Subsection Stacks of polynomial growth
\begin{subsection}{Stacks of polynomial growth}

The aim of this section is to try to fit the considerations of \cite{behrend07}
into our setup. Hence, if $\sX$ is an algebraic stack locally of finite type
(and as usual with affine stabilisers), a \Definition{stratification} of $\sX$
will mean a countable set $(\sX_i)$ of locally closed substacks of finite type
such that $\sX$ is the disjoint union of them and $\dim(\sX_i) \to -\infty$ as
$i \to \infty$. (Note that we do \emph{not} assume for instance that the closure
of a stratum is a union of strata.) Following \cite{behrend07} we say that $\sX$
is \Definition{essentially of finite type} if it has a stratification. The
stratification is \Definition{uniformly of polynomial growth} if the sum
$\sum_i\overline{\{\sX_i\}}$ (we shall from now one dispense with the completion
symbol) converges in $\Kspol{\k}$. In concrete terms this means that there is a
polynomial $P(x)$ such that $\uwt_n(\sum_{i\in S}\sX_i) \le P(|n|)$ for all $n$
and all finite subsets $S$ of the index set as convergence in the filtration
topology is assured by (\ref{completion continuity}). (This is rather the
summability than the convergence series version but they are easily seen to be
equivalent.) We say that $\sX$ is of \Definition{polynomial growth} if it has
some stratification $(\sX_i)$ which is uniformly of polynomial growth. We then
want to define $\{\sX\} \in \Kspol{\k}$ as $\sum\{\sX_i\}$. This is indeed
independent of the choice of stratification. To prove that, it is enough to show
equality in $\Kschemecompl{\k}$ (as the map $\Kspol{\k} \to \Kschemecompl{\k}$
is injective). Hence, if we simply assume that $(\sX_i)$ and $(\sY_j)$ are
exhaustions of $\sX$ we see that both sums $\sum_i\{\sX_i\}$ and
$\sum_j\{\sY_j\}$ are equal to $\sum_{i,j}\{\sX_i\cap\sY_j\}$ as for each $i$
the stratification $(\sX_i\cap\sY_j)$ of $\sX_i$ is a finite stratification as
$\sX_i$ is of finite type (and similarly for the stratification
$(\sY_j\cap\sX_i)$ of $\sY_j$). This gives $\{\sX_i\}=\sum_j\{\sX_i\cap\sY_j\}$
and $\{\sY_j\}=\sum_i\{\sY_j\cap\sX_i\}$. We then have a few basic properties.
%;;%%%Proposition
\begin{proposition}
\part[i] Suppose that $\sX$ is an algebraic $\k$-stack locally of finite type
and $\sY$ is a closed substack with complement $\sU$. If the stacks $\sY$ and
$\sU$ are of polynomial growth then so is $\sX$ and we have
\begin{displaymath}
\{\sX\} = \{\sY\} + \{\sU\} \in \Kspol{\k}.
\end{displaymath}

\part[ii] If $\sX$ is a stack of polynomial growth over a finite field $\k$, then
$\Tr(\chi_c\{\sX\})=\cnt(\sX)$.
\begin{proof}
Suppose that $\sY$ resp.\ $\sU$ have stratifications $(\sY_i)$ resp.\ $(\sU_i)$
of uniform polynomial growth. Their union then gives a stratification of $\sX$
of uniform polynomial growth. This proves \DHrefpart{i}.

As for \DHrefpart{ii} we have, by definition of the point count, that
$\cnt(\sX)=\sum_i\cnt(\sX_i)$. Hence by Proposition \ref{Finite type point
counting} we have
\begin{displaymath}
\Tr(\{\sX\}) = \sum_i\Tr(\{\sX_i\})=\sum_i\cnt(\sX_i) = \cnt(\sX).
\end{displaymath}
\end{proof}
\end{proposition}
\begin{remark}
\part It would seem reasonable to demand a version for the first part that
involves an infinite decomposition of $\sX$ into disjoint substacks. However,
one would to begin with have to impose some uniformity on the dimensions of the
strata of the stratification of the individual pieces. A second problem is that
if $(\sX_i)$ is a stratification of $\sX$ then it is not clear that one can
bound $\sum_i\uwt(\{\sX_i\})$ in terms of $\uwt_n(\{\sX\})$ which would be
needed to bound $\sum_i\uwt_n(\{\sX_i\})$ in terms of
$\sum_{i,j}\uwt_n(\{\sX_{i,j}\})$, where $(\sX_{i,j})$ is a stratification of
$\sX_i$ of uniform polynomial growth (and $(\sX_i)$ is a decomposition of
$\sX$).

\part Another problem with our definition is that (even open) substacks of a
stack of polynomial growth may not be of polynomial growth. This can be
mitigated by modifying the notion of convergence. Unfortunately under that
modification the Euler characteristic is no longer necessarily of polynomial
growth.
\end{remark}
We finish this subsection by verifying that the moduli stacks of
\cite{behrend07}, stacks of bundles over a smooth and proper curve for a split
semi-simple group $G$.
%;;%%%Proposition
\begin{proposition}
Let $G$ be a split semi-simple $\k$-group scheme and $C$ a smooth and proper
curve over $C$. Then the algebraic stack $\Bun_{G,C}$ of $G$-torsors over $C$ is
of polynomial growth.
\begin{proof}
This is just going into the details of \cite[Lemma 5.5]{behrend07} to see that
the convergence is actually polynomial using \cite[\S8]{behrend91::lefsc},
particularly \cite[\S8.4]{behrend91::lefsc} to get the requisite
information. For the reader's (and the author's) convenience we shall only treat
the notationally (and conceptually) simpler case of $G=\SL_n$. Hence we stratify
$\Bun_{\SL_n,C}$ by the strata of fixed ranks and degrees of the pieces of the
Harder-Narasimhan filtration. Hence they are given by one sequence
$(d_1,d_2,\dots,d_k)$ of positive integers giving the ranks and one increasing
sequence $(m_1,\dots,m_k)$ of integers representing the degrees of the pieces of
the filtration. That the whole vector bundle comes from an $\SL_n$-bundle means
that $\sum_id_i=n$ and that we have a Harder-Narasimhan filtration means that
$m_1/d_1 > m_2/d_2 > \dots > m_k/d_k$. There is a small problem in that on a
given stratum the Harder-Narasimhan filtration may not exist. However, the map
from the algebraic stack of vector bundles with a Harder-Narasimhan filtration
of the given type to the stratum is a universal homeomorphism. Generally, if
$\sX \to \sY$ is a universal homeomorphism between algebraic stacks of finite
type we have that $\chi_c(\{\sX\})=\chi_c(\{\sY\})$ in $\Kcpol{\k}$, a fact
which follows immediately by reduction to the spatial case. Hence, we may and
shall ignore this problem. Now, $\Bun_C^{d,m}$ is the stack of flags $0\subset
\sE_1 \subset\dots\subset\sE_k$ of vector bundles on $C$ with $\sE_i/\sE_{i-1}$
semi-stable of rank $d_i$ and degree $m_i$ (and where we have $\sum_id_i=n$ and
$m_1/d_1 > m_2/d_2 > \dots > m_k/d_k$). We now consider the map $\Bun_C^{d,m}
\to \Bun_C^{d',m'}\times\Bun_C^{d_k,m_k}$ mapping a flag to $(0\subset \sE_1
\subset\dots\subset\sE_{k-1},\sE_k/\sE_{k-1})$ and $d'=(d_1,\dots,d_{k-1})$ and
$m'=(m_1,\dots,m_{k-1})$. What is lost under this map is the extension of
$\sE_{k-1}$ by $\sE_k/\sE_{k-1}$ given by $\sE_k$. This makes $\Bun_C^{d,m}
\to \Bun_C^{d',m'}\times\Bun_C^{d_k,m_k}$ a Picard stack under the Baer sum of
extensions which is described by the perfect complex
$R\pi_*\sHom(\sF,\sE_{k-1})[1]$, where $0\subset \sE_1
\subset\dots\subset\sE_{k-1}$ is the universal flag over $\Bun_C^{d',m'}\times C$,
$\sF$ the universal vector bundle over $\Bun_C^{d_k,m_k}\times C$ and $\pi$ is
the projection $\Bun_C^{d',m'}\times\Bun_C^{d_k,m_k}\times C \to
\Bun_C^{d',m'}\times\Bun_C^{d_k,m_k}$. We may then stratify the base so that
$\pi_*\sHom(\sF,\sE_{k-1})$ is free of rank $a$, say, and
$R^1\pi_*\sHom(\sF,\sE_{k-1})$ is free of rank $b$, say, and so that
$R\pi_*\sHom(\sF,\sE_{k-1})$ is the sum of its cohomology. This gives that over
such a stratum $\sX$, $\pi$ is isomorphic to $\Bclass
\add^a\times\add^b\times$. Using additivity of $\{-\}$ we get 
\begin{displaymath}
\{\Bun_C^{d,m}\}=\Lclass^{b-a}\{\Bun_C^{d',m'}\times\Bun_C^{d_k,m_k}\}=
\Lclass^{b-a}\{\Bun_C^{d',m'}\}\{\Bun_C^{d_k,m_k}\}
\end{displaymath}
and $a-b$ is simply the (relative) Euler characteristic of
$\sHom(\sF,\sE_{k-1})$ which equals 
\begin{displaymath}
\sum_{i<k}d_km_i-d_im_k
\end{displaymath}
by the Riemann-Roch formula. Iterating we get
\begin{displaymath}
\{\Bun_C^{d,m}\}=\Lclass^{-\chi}\{\Bun_C^{d_1,m_1}\}\cdots\{\Bun_C^{d_k,m_k}\},
\end{displaymath}
where
\begin{displaymath}
\chi=\sum_{i<j}d_id_j\left(\frac{m_j}{d_j}-\frac{m_i}{d_i}\right)+d_id_j(1-g(C)).
\end{displaymath}
Now, for fixed $d$ $\sum_{i<j}d_id_j(\frac{m_j}{d_j}-\frac{m_i}{d_i})$ is a
linear function whose kernel clearly meets the cone given by the conditions
$m_1/d_1 > m_2/d_2 \ge \dots \ge m_k/d_k$ and $\sum_im_i=0$ only in $0$. Hence
the number of $m$'s which give a fixed value of $\chi$ is polynomially bounded
and only takes on a finite number of negative values. As also the isomorphism
class of each $\Bun_C^{d_i,m_i}$ only depends on the residue of $m_i$ modulo
$d_i$ (provided that the field is large enough so that $C$ has a line bundle of
degree $1$, field extensions are allowed however) we get that sum
$\sum_{d,m}\{\Bun_C^{d,m}\}$ is uniformly polynomially convergent.
\end{proof}
\end{proposition}
\begin{remark}
The cohomology of the moduli spaces in question ultimately comes from the
cohomology of the group in question, its loop space and its classifying space
(cf., \cite{atiyah83::yang+mills+rieman}). The polynomial growth condition is
thus related to the polynomial growth of the cohomology of the classifying space
and loop space of the group. Note that there are finite CW-complexes for which the
cohomology of its loop space grows exponentially. This may or may not be an
argument against seeing the polynomial growth condition as a natural condition.
\end{remark}
\end{subsection}
%;;%%Subsection The real $\chi_c$ of algebraic stacks
\begin{subsection}{The real $\chi_c$ of algebraic stacks}

Note that so far we have been defining $\chi_c$ for an algebraic stack in
terms of $\chi_c$ of algebraic spaces. Using the theory of cohomology of compact
support for algebraic $\k$-stacks of finite type introduced in
\cite{laszlo08::artin+ii} we shall now see that $\chi_c$ may indeed be
interpreted as an actual (convergent) Euler characteristic of the stack. We
start by checking that some expected properties are indeed true. Following
\cite{laszlo08::artin+ii} we shall only consider $\ell$-adic cohomology
(although it seems reasonably clear that their arguments can be carried through
also in the case of mixed Hodge structures). We first treat the case of an
algebraic stack of finite type (we shall use $H^i_c(\sX)$ to denote the
cohomology of the extension of $\sX$ to a separable closure of $\k$ with
$\Q_\ell$-coefficients).
%;;%%%Proposition finite type basics (ii)Dimension bound
\begin{proposition}\label{finite type basics} 
Let $\sX$ be an algebraic stack of finite type.

\part[i] If $\sU \subseteq \sX$ is an open substack and $\sZ$ its complement, then
we have a long exact sequence
\begin{displaymath}
\dots \to H^i_c(\sU) \to  H^i_c(\sX) \to H^i_c(\sZ) \to H^{i+1}_c(\sU) \to\dots
\end{displaymath}

\part[ii]\label{Dimension bound} If $\dim\sX \le N$, then $H^i_c(X)=0$ for $i > 2N$.

\part[iii] $H^i_c(\sX)$ is mixed of weight $\le i$.

\part[iv] $\dim H^i_c(\sX)$ has polynomial growth.
\begin{proof}
\DHrefpart{i} follows from \cite[4.9.0.1]{laszlo08::artin+ii} and a passage to
the limit. As for \DHrefpart{ii} we can use \DHrefpart{i} to reduce to the case
when $\sX$ is smooth and then the result follows from duality. Further, to prove
\DHrefpart{iii} we can use \DHrefpart{i}, and the fact that mixed sheaves are
closed under subobjects, quotients and extensions, to reduce to the case when
$\sX$ is a global quotient. In that case we can find a vector bundle over $\sX$
with an open substack which is spatial and with complement of
arbitrarily high codimension. As the cohomology of a vector bundle is the
cohomology of the base shifted the same amount in degree and weight we can use
\DHrefpart{ii} and \DHrefpart{i} to reduce to the case an algebraic space in
which case it is Deligne's fundamental result, \cite{deligne80::la+weil}.

Finally, to prove \DHrefpart{iv} we are reduced to the case when $\sX$ is smooth
and a global quotient $[X/\GL_n]$. We are then by duality reduced to proving
that $\dim H^i(\sX)$ has polynomial growth which follows from the spectral sequence
\begin{displaymath}
H^i(\Bclass \GL_n)\Tensor H^j(X) \implies H^{i+j}(\sX).
\end{displaymath}
\end{proof}
\end{proposition}
Now, \cite{laszlo08::artin+i} defines cohomology with compact support for finite
coefficients for algebraic stacks \emph{locally} of finite type. It is easy to
see that this cohomology (for constant coefficients say) is the inductive limit
of the cohomology over all open substacks of finite type. This property is
however destroyed in general when one passes to $\ell$-adic or
$\Q_\ell$-coefficients. It therefore seems better to \emph{define} $H^i_c(\sX)$
for a stack only locally of finite type as the inductive limit over open
substacks of finite type. The only thing to decide is in which sense the limit
is to be taken, we choose to simply to take the ind-system defined by the open
substacks of finite type. This makes sense both in the Galois and the mixed
Hodge structure case. (This whole discussion is somewhat moot as in the cases
that we shall consider each $H^i_c(\sX)$ will be an eventually constant system,
where an ind-system $\{\sH_\alpha\}$ is \Definition{eventually constant} if
there is an index $\alpha$ so that for any morphism $\alpha \to \beta$, the
induced map $\sH_\alpha \to \sH_\beta$ is an isomorphism. This is turn means
that the system is isomorphic as an ind-object to a constant system.)
%;;%%%Theorem Real Euler characteristic
\begin{theorem}\label{Real Euler characteristic}
\part[i] If $\sX$ is essentially of finite type, then for each $i$, $H^i_c(\sX)$ is
eventually constant.

\part[ii] If $\sX$ is of uniform polynomial growth and $\dim H^i_c(\sX)$ grows
polynomially, then
\begin{displaymath}
\chi_c(\{\sX\}) = \sum_i(-1)^i\{H^i_c(\sX)\} \in \Kcpol{\k}.
\end{displaymath}
\begin{proof}
Let $\{\sX_i\}$ be a stratification of $\sX$ and $\{\sU_\alpha\}$ the family of
open substacks of finite type. Given $N$ there is only a finite number of
$\sX_i$ of dimension $>N$. Each of them can be covered by a finite number of the
$\sU_\alpha$ and as $\{\sU_\alpha\}$ is inductive there is a $\sU_\alpha$ whose
complement has dimension $\le N$ and hence by Proposition \ref{finite type
basics} for any $\sU_\alpha \subseteq \sU_\beta$, $H^i_c(\sU_\alpha) \to
H^i_c(\sU_\beta)$ is an isomorphism which proves \DHrefpart{i}.

Continuing, the assumption on polynomial growth and the fact that $H^i_c(\sX)$
is of weight $\le i$ implies that $\sum_i(-1)^i\{H^i_c(\sX)\}$ converges and to
show that it is equal to $\chi_c(\{\sX\})$ it is enough to show this in
$\Kcompl{\k}$. This amounts to showing that the weight $n$-parts are equal. Now,
using the notation of \DHrefpart{i}, the weight $n$ part of
$\sum_i(-1)^i\{H^i_c(\sX)\}$ is equal to the weight $n$-part of
$\sum_i(-1)^i\{H^i_c(\sU_\alpha)\}$ with $N<<0$. On the other hand, the weight
$n$-part of $\chi_c(\{\sX\})$ is equal to the weight $n$ part of 
\begin{displaymath}
\chi_c(\sum_i\{\sX_i\cap\sU_\alpha\}) = \chi_c(\sU_\alpha)
\end{displaymath}
as all but a finite number of the $\sX_i\cap\sU_\alpha$ are empty. Hence we are
reduced to the case when $\sX$ is of finite type. Cutting $\sX$ into pieces we
may assume that $\sX$ is a global quotient. This means that there are vector
bundles over $\sX$ with non-spatial locus of arbitrarily high codimension. This
high codimension locus may be ignored by (\ref{completion continuity}) and
(\ref{Dimension bound}). Hence we are reduced to the case when $\sX$ is a space
in which case it clear (by definition).
\end{proof}
\end{theorem}
\end{subsection}
\end{section}
%;;%Section The motivic class and the refined Euler characteristic
\begin{section}{The motivic class and the refined Euler characteristic}

When the base field has characteristic zero we can use the presentation of
$\Kscheme{\k}$ (and consequently of $\Kstack{\k}$) obtained by Bittner (cf.,
\cite{bittner04::euler}) to get more refined cohomological invariants of classes
in $\Kscheme{\k}$. We start by defining the recipients for these invariants. As
we are going to use the Bittner representation \emph{in this section $\k$ will
be a field of characteristic zero}. For a general field $\k$ (of characteristic
zero) we let $\Gcoh{\k}$ be the group generated by isomorphism classes in the
category of direct limits of bounded constructible $\Z_p$-complexes (in the
derived category) over finite type ``thickenings'' of $\k$ (where $p$ is a fixed
prime). We introduce the relations $\{A\Dsum B\}=\{A\}+\{B\}$ but \emph{no}
relations for distinguished triangles. In case $\k=\C$ we alternatively let
$\Gcoh{\k}$ be generated by isomorphism classes of bounded $\Z$-complexes $A$
with finitely generated cohomology together with a (polarisable) Hodge structure
of weight $k$ on $H^k(A)$ and the same relations as before (and when $\k=\R$ the
same plus involutions as above). We define a ring structure on $\Gcoh{\k}$ by
putting $\{A\}\{B\}=\{A\Tensor^L_{\Z_p}B\}$ (resp., of course,
$\{A\Tensor^L_{\Z}B\}$). There is an increasing \Definition{degree filtration}
on $\Gcoh{\k}$ with $\Fil^k\Gcoh{\k}$ being generated by complexes concentrated
in degrees $\le k$ and we let $\Gcompl{\k}$ be the completion of $\Gcoh{\k}$
with respect to this filtration. As the filtration is multiplicative (i.e.,
$\Fil^k\cdot\Fil^\ell\subseteq \Fil^{k+\ell}$) we have a ring homomorphism
$\Gcoh{\k}\to\Gcompl{\k}$. Before our first major result we need the following
lemma.
%;;%%Lemma Filtration properly spanned
\begin{lemma}\label{Filtration properly spanned}
$\Fil^n\Kscheme{\k}$ is spanned by classes of smooth and proper varieties of
dimension $\le n$.
\begin{proof}
We prove this by induction over $n$ with $n=-1$ being obvious. Assume now that
$\dim U\le n$, we then want to express it as a linear combination of smooth and
proper varieties of dimension $\le n$. Choosing first a compactification of $U$
and then a resolution of singularities expresses the class as a linear
combination of an $n$-dimensional smooth and proper variety and varieties of
smaller dimension which finishes the proof by induction.
\end{proof}
\end{lemma}
%;;%%Proposition Refined Euler
\begin{proposition}\label{Refined Euler}
\part[i] There is a ring homomorphism \map{\chi_G}{\Kscheme{\k}}{\Gcoh{\k}}
characterised by $\chi_G(\{X\})=\{R\Gamma(X)\}$ where $X$ is smooth and proper.

\part[ii] $\chi_G$ extends (uniquely) to a continuous ring homomorphism
\map{\chi_G}{\Kschemecompl{\k}}{\Gcompl{\k}}. In particular it extends
(uniquely) to a ring homomorphism \map{\chi_G}{\Kstack{\k}}{\Gcompl{\k}}.

\part For each $k$ there is a group homomorphism
\map{H^k(-)\Tensor\Q}{\Gcompl{\k}}{\Kcompl{\k}} characterised by continuity and the property
$H^k(\{C\})\Tensor\Q=\{H^k(C)\Tensor\Q\}$. We have for $x \in \Kschemecompl{\k}$
that $H^k(x)\Tensor\Q$ is pure of weight $k$ and
$\chi(x)=\sum_i(-1)^kH^k(\chi_G(x))$. In particular $H^k(\chi_G(x))\Tensor\Q$ is the
weight $k$ part of $(-1)^k\chi(x)$.
\begin{proof}
By \cite[thm.\ 3.1]{bittner04::euler} $\Kscheme{\k}$ is generated by classes
$\{X\}$, $X$ smooth and proper, and their relations are generated by
$\{\mathrm{Bl}_YX\}-\{E\}=\{X\}-\{Y\}$, where $X$ is smooth and proper,
$Y\subseteq X$ a smooth and proper subvariety, $\mathrm{Bl}_YX$ is the blowing up
of $X$ along $Y$. Now, there is the well-known isomorphism of complexes
$R\Gamma(\mathrm{Bl}_YX)\Dsum R\Gamma(Y)\iso R\Gamma(X)\Dsum R\Gamma(Y)$ (e.g.,
$R\Gamma(-)$ factors through the category of integral motives and the
corresponding isomorphism of motives is \cite[Cor.\ \S7,Cor.\
\S9]{manin68::corres}). This proves \DHrefpart{i}.

For \DHrefpart{ii} we first note that $\chi_G$ takes $\Lclass$ to an invertible
element in $\Gcompl{\k}$. Indeed, $\Lclass=\{\P^1\}-1$ and
$\chi_G(\{\P^1\})=1+\{\Z(-1)[-2]\}$ (resp.\ $\Z_p(-1)[-2]$) and $\{\Z(-1)[-2]\}$
is invertible, with inverse $\{\Z(1)[2]\}$ already in $\Gcoh{\k}$. This gives us
a ring homomorphism $\Kscheme{\k}[\Lclass^{-1}]\to \Gcompl{\k}$ and we finish if
we can show that it is continuous in the filtration topologies. This however
follows directly from Lemma \ref{Filtration properly spanned} and the facts that
$\chi_G(\{X\}/\Lclass^m)=\{R\Gamma(X)(m)[2m]\}$ and that $R\Gamma(X)$ is of
amplitude $[0,2n]$ if $X$ is $n$-dimensional. The only thing left is then the
uniqueness of the extension of $\chi_G$ to $\Kstack{\k}$ but it follows from the
fact that $\Kscheme{\k} \to \Kstack{\k}$ is a localisation.

For the last part, as $H^k(-)$ is additive on complexes we get an induced map
\map{H^k(-)\Tensor\Q}{\Gcoh{\k}}{\Kcoh{\k}} which clearly is continuous and
hence extends to $\Gcompl{\k}$. Now both $H^k(\chi_G(-))$ is additive and
continuous so to verify that it maps to the weight $k$ part of $\Gcompl{\k}$ it
is enough to very it on classes $H^k(\chi_G(\{X\}))\Tensor\Q$ where $X$ is smooth
and proper. In that case $H^k(\chi_G(\{X\}))\Tensor\Q=\{H^k(X,\Q)\}$ (resp.\
with $\Q_p$-coefficients) and $H^k(X,\Q)$ is indeed pure of weight $k$. This
implies that the sum $\sum_i(-1)^kH^k(\chi_G(x))$ converges and gives a
continuous additive function on $\Kschemecompl{\k}$ and so does $\chi$. Hence to
verify that they are equal it suffices to do it on $\{X\}$, $X$ smooth and
proper, in which case it is clear.
\end{proof}
\end{proposition}
\begin{remark}
Instead of the Bittner presentation one could use earlier results (cf.,
\cite{gillet96::descen+k}) that associates a virtual motive to a class of
$\Kcoh{\k}$.
\end{remark}
In order to detect non-trivial elements of $\Gcoh{\k}$ (or $\Gcompl{\k}$) it is
useful to pass to the underlying abelian groups. Hence, we let $\Gcplx$ be the
group generated by isomorphism classes of bounded complexes (in the derived
category) of finitely generated $\Z$-complexes (resp., of course,
$\Z_p$-complexes) by the same relations as before and $\Gcmpl$ the corresponding
degree completion. Again they have a ring structure given by the derived tensor
product. We also let $\Gab$ be the corresponding group based on modules rather
than complexes but we do not give it a multiplicative structure. We have a group
homomorphism $\Gab \to \Gcplx$ given by $\{M\} \mapsto \{M[0]\}$ and
\map{H^k}{\Gcmpl}{\Gab} given by $\{C\}\mapsto \{H^k(C)\}$. We continue by
defining $1:=\{\Z\} \in \Gab$ (which maps to the identity element of $\Gcplx$)
and $\alpha_{p,n}:=\{\Z/q^n\}$, where $n>0$ and $q$ is a prime. Finally, for $M$
an abelian group we let $M[t,t^{-1}]$ be the group of finite formal sums
$\sum_{i\in\Z}m_it^i$, $m_i \in M$ and $M((t^{-1}))$ the group of formal sums
$\sum_{i\in \Z}m_it^i$ which possibly have an infinite number of negative power
terms.
%;;%%Proposition Complex classification
\begin{proposition}\label{Complex classification}
\part $\Gab$ is the free abelian group on $1$ and the $\alpha_{q,n}$, where $q$
runs over all primes (resp.\ is equal to $p$).

\part The maps \map{H^\bullet}{\Gcplx}{\Gab[][t,t^{-1}]} and
\map{H^\bullet}{\Gcmpl}{\Gab((t^{-1}))} given by
$H^\bullet(\{C\})=\sum_iH^i(\{C\})t^i$ are isomorphisms.

\part The multiplicative structure on $\Gab[][t,t^{-1}]$ and $\Gab((t))$ induced
by these isomorphisms is characterised by the conditions that the maps
$\Z[t,t^{-1}]\to \Gab[t,t^{-1}]$ and $\Z((t^{-1}))\to \Gab((t))$ induced by $m
\mapsto m\cdot 1$, $t \mapsto t$ and the second map is continuous are ring
homomorphisms and
\begin{displaymath}
\alpha_{q,m}\alpha_{q',n}=
\begin{cases}
0&\text{if }q\ne q'\\
\alpha_{q,\min(m,n)}(1+t^{-1})&\text{if }q= q' 
\end{cases}
\end{displaymath}
together with the fact that the multiplication is continuous in the $t^{-1}$-adic
topology in the $((t^{-1}))$-case.
\begin{proof}
As the rings $\Z$ and $\Z_p$ have global dimension $1$ the complexes in question
are isomorphic to the direct sum of their (shifted) cohomology modules. This
gives the second part. The first part follows from the classification of
finitely generated modules over the rings in question and the last from the
basic calculation of tensor products and $\mathrm{Tor}$-groups.
\end{proof}
\end{proposition}
\begin{remark}
\part The conclusion should be of course thought of as saying that the usual
rank and torsion invariants of the cohomology of a smooth and proper variety are
in fact invariants of its class in $\Kscheme{\k}$ (and in fact of its class in
$\Kschemecompl{\k}$). The precise structure of the groups are being mentioned
mainly to emphasise how very much larger they are than the ordinary Grothendieck
groups (which are just isomorphic to $\Z$).

\part If one works rationally instead of integrally, then the refined Euler
characteristic can be recovered from the ordinary Euler class in
$\Kcoh{\k}$. Indeed, $R\Gamma(X)$ is a pure complex so is isomorphic to the sum
of its (shifted) cohomology and then we can use Theorem \ref{Refined
Euler}. Hence, the extra information provided by the refined Euler
characteristic is torsion and integrality information.
\end{remark}
There are a number of variants of this technique some of which will appear in
the next result. To prepare for it let us note that there is an obvious additive
Bittner type presentation of $\Kscheme{\k}[\Lclass^{-1}]$: It is generated by
classes $\{X\}/\Lclass^n$, $n\ge 0$, and relations
$\{\mathrm{Bl}_YX\}/\Lclass^n-\{X\}/\Lclass^n=\{E\}/\Lclass^n-\{Y\}/\Lclass^n$
as before and $\{\P^1\times
X\}/\Lclass^{n+1}-\{X\}/\Lclass^{n+1}=\{X\}/\Lclass^n$ for smooth and proper $Y
\subseteq X$. Also we shall put $\widehat\Z:=\ili_n\Z/n=\prod_p\Z_p$ and
$H^*(X,\widehat\Z):=\prod_pH^*(X,\Z_p)$. 
%;;%%Theorem 
\begin{theorem}
\part[i] Let $A_\k$ be the group generated by isomorphism classes of algebraic
group schemes $A$ over $\k$ whose connected component is an abelian variety and
whose group of geometric connected components has a finitely generated group
of geometric points and with relations $\{A\Dsum B\}=\{A\}+\{B\}$. Then there is
a (unique) group homomorphism \map{\Pic_\k}{\Kscheme{\k}}{A_{\k}} for which
$\Pic_\k(\{X\})=\{\Pic(X)\}$ for a smooth and proper $X$.

\part[ii] There is an extension of $\Pic_\C$ to $\Kschemecompl{\C}$ vanishing on
$\Fil^0\Kschemecompl{\C}$ (which we shall also call $\Pic_\C$).

\part[iii] Let $\Gab[\k]$ be the group generated by isomorphism classes of
\'etale algebraic group schemes $A$ whose group of geometric connected
components has a finitely generated group of geometric points and with relations
$\{A\Dsum B\}=\{A\}+\{B\}$. Then there are continuous group homomorphisms
\map{\NS^k}{\Kschemecompl{\k}}{\Gab[\k]} characterised by
$\NS^k(\{X\}/\Lclass^n)=\NS^{k+n}(\{X\})$ and $\NS^k(\{X\})=\{\NS^k(X)\}$, for
$X$ smooth and proper, where $\NS^k(X)(\bk)$ is the group of algebraic
codimension $k$-cycles on $X_\bk$ modulo homological equivalence with its
natural $\Gal(\bk/\k)$-action.
\begin{proof}
Starting with \DHrefpart{i} it follows directly from the Bittner presentation
and the computation of the Picard scheme of a blow up. As for \DHrefpart{ii} we
define $\Pic'(\{X\}/\Lclass^n)$, $X$ smooth and proper, as follows: We let it be
the class of $A^{0,n}_X\Dsum A^{c,n}_X$, where $A^{c,n}_X$ is the inverse image
of the classes of type $(n+1,n+1)$ in $H^{2n+2}(X,\C)$ under the map
$H^{2n+2}(X,\Z) \to H^{2n+2}(X,\C)$. Similarly, we let $A^{0,n}_X$ be the Weil
intermediate Jacobian (cf., \cite{weil52::picar}) associated to the
(polarisable) Hodge structure on $H^{2n+1}(X,\Z)$. Now $\Pic(X)$, as any complex
commutative algebraic group whose group of components is an abelian variety, is
the direct product of its connected component and its group of components (as we
are working over an algebraically closed field). Furthermore, $A^{0,1}_{X}$ is
the Jacobian, i.e., the connected component of $\Pic(X)$, and $A^{c,1}(X)$ is
the group of components of $\Pic(X)$ by Lefschetz theorem. It is furthermore
clear that $\{\P^1\times X\}/\Lclass^{n+1}-\{X\}/\Lclass^{n+1}$ is mapped to the
same element as $\{X\}/\Lclass^n$ which shows that we $\Pic'$ is given as a map
from $\Kscheme{\k}[\Lclass^{-1}]$. Finally, as $R\Gamma(X)$ is of amplitude
$[0,2\dim X]$ it is clear (using Lemma \ref{Filtration properly spanned})that
$\Fil^0\Kscheme{\k}$ is mapped to zero which in particular gives us an extension
to $\Kschemecompl{\k}$ with the same property.

Turning to \DHrefpart{iii}, by definition $\NS^k(X)$, for $X$ smooth and
proper, is the subgroup of $H^{2k}(X,\widehat\Z(k))$ spanned by the cycle
classes of $k$-codimensional subvarieties of $X_\bk$. It is clearly stable by
the Galois action and hence gives rise to an \'etale group scheme. The next step
is to show that it is finitely generated. If not there is a countable set
$\{Y_i\}$ of $k$-codimensional subvarieties whose cycle classes span a subgroup $M$
of $H^{2k}(X,\widehat\Z(k))$ with the property that there is no finite subset of
$Y_i$ whose cycle classes span $M$. The subvarieties $Y_i$ may be defined over a
common countably generated algebraically closed subfield $K$ of $\bk$ which then
may be embedded in $\C$. We may then replace $\bk$ by $\C$ in which case the
cycle class map factors through $H^{2k}(X,\Z) \to H^{2k}(X,\widehat\Z(k))$ which
shows that $M$ is finitely generated which is a contradiction.

What remains to be verified are the Bittner relations and $\{\P^1\times
X\}/\Lclass^{n+1}-\{X\}/\Lclass^{n+1}=\{X\}/\Lclass^n$. These follow immediately
from the formulas for the cohomology and Chow groups of blowing ups and products
with $\P^1$ and the compatibility of these formulas with the cycle map.
\end{proof}
\end{theorem}
\begin{remark}
\part Instead of the Weil intermediate Jacobian we could have used the maximal
abelian subvariety of the usual (=''Griffiths'') intermediate Jacobian
$H^{2n+1}(X,\C)/H^{2n+1}(X,\Z)+F^nH^{2n+1}(X,\C)$.

\part The maps $\NS^k$ will be used in \cite{ekedahl08::stack} to fit a
counterexample of Swan to Noether's problem into our context.
\end{remark}
\begin{corollary}
The rings $\Kscheme{\k}$, $\Kscheme{\k}[\Lclass^{-1}]$, $\Kstack{\k}$
and $\Kschemecompl{\k}$ contain zero-divisors.
\begin{proof}
Let us recall the setup of \cite{poonen02::groth} (where Poonen proves the
statement for $\Kscheme{\k}$). We have an abelian variety $A$ over $\k$ with
endomorphism ring (over $\k$ as well as geometrically) equal to the ring of
integers $R$ in a number field. Furthermore, we have a non-principal ideal $I$
in $R$ and put $B:=I\Tensor_RA$ (which in particular has the property that
$\Hom(A,B)=I$ as $R$-module). We have that $B^n\iso I^{\dsum n}\Tensor_RA\iso
A^{n-1}\Dsum A\Tensor_RI^n$ and if we choose $n$ such that $I^n$ is principal
($n=2$ in Poonen's specific example) we get $B^n\iso A^n$ which gives
$0=(\{A\}^{n-1}+\{A\}^{n-2}\{B\}\cdots+B^{n-1})(\{A\}-\{B\})$ in $\Kscheme{\k}$
and we want to show that both factors have non-zero images in any of the rings
in question and as all of them map compatibly to $\Kschemecompl{\k}$ it is
enough to show that their images in that ring are non-zero. For the first factor
this is trivial, the weight $2\dim A(n-1)$-part of its Euler characteristic is
non-trivial. We are thus left with the problem on whether $\{A\}=\{B\} \in
\Kschemecompl{\k}$. Now, to show that they are different it is enough to show
that they are distinct modulo $\Fil^0$ and hence as elements of
$\Kscheme{\k}[\Lclass^{-1}]$. For this we may by standard limit arguments assume
that $\k=\C$. As $\Pic(\{A\})=\{A\} \in A_\C$ and the same for $B$ we are
reduced to showing that $A\Dsum C\iso B\Dsum C$ for some abelian group scheme
$C$ whose connected component is an abelian variety. Such an isomorphism induces
an isomorphism $\Hom(A,A)\Dsum\Hom(A,C)\iso\Hom(A,B)\Dsum\Hom(A,C)$ of finitely
generated $R$-modules. Putting $M:=\Hom(A,C)$ this gives an isomorphism $R\Dsum
M \iso I\Dsum M$ but we have cancellation for finitely generated $R$-modules so
that we get an isomorphism $R \iso I$ which contradicts the assumption that $I$
is non-principal.
\end{proof}
\end{corollary}
It is natural to ask in what generality the last argument works. The general
question is to what extent a group scheme is determined by its class in
$A_\k$. This is the same problem as asking whether $A\Dsum C \iso B\Dsum C$
implies that $A\iso B$. We restrict ourselves to the case when $A$ and $B$ are
abelian varieties and by restricting to identity components we may assume that
$C$ is also an abelian variety.
%;;%%Proposition
\begin{proposition}
Assume that $A$, $B$ and $C$ are abelian varieties and $A\Dsum C\iso
B\Dsum C$. Then there is a locally free rank $1$ module $I$ over $R:=\End(A)$
such that $B$ is isomorphic to $I\Tensor_RA$ so that in particular $A\iso B$
precisely when $R\iso I$ and a torsion free finitely
generated $R$-module $M$ such that $R\Dsum M\iso I\Dsum M$. In particular if $R$
is hereditary and no simple factor of $R\Tensor \Q$ is a totally definite
quaternion algebra over its centre then $A\iso B$.
\begin{proof}
We may assume that $\k$ is finitely generated. We get for any prime $p$ an
isomorphism $T_pA\Dsum T_pC \iso T_pB\Dsum T_pC$ of modules over the Galois
group $\sG$ of $\k$. The Krull-Schmidt theorem then implies that we have an
isomorphism $T_pA \iso T_pB$ of $\sG$-modules. Let $I:=\Hom(A,B)$. By the Tate
conjecture \cite[Thm 4]{faltings83::endlic+variet+zahlk} we have
$I_p:=I\Tensor\Z_p=\Hom_\sG(T_pA,T_pB)$ and in particular, using the isomorphism
$T_pA \iso T_pB$ we get that $I_p \iso R_p:=R\Tensor\Z_p$ as $R_p$-modules so
that $I$ is a locally free $R$-module of rank $1$. This implies first that the
evaluation map $\Hom_\sG(T_pA,T_pB)\Tensor_{R_p}T_pA \to T_pB$ and second that
$T_p(I\Tensor_RA)=I\Tensor_RT_p(A)$. Together these facts imply  that the
tautological map $I\Tensor_RA \to B$ induces an isomorphism when $T_p(-)$ for
all $p$ is applied and hence is an isomorphism which gives the first
part. Putting $M:=\Hom(A,C)$ we get an isomorphism of $R$-modules $R\Dsum
M=\Hom(A,A\Dsum C)\iso\Hom(A,B\Dsum C)=I\Dsum M$.

Finally, when $R$ is hereditary the fact that $M$ is torsion free implies that
$M$ is projective. The condition on $R\Tensor\Q$ implies that it fulfils the
Eichler condition relative to $\Z$ (see \scite[Def.\ 38.1]{Re75}) and hence we
have cancellation for projective modules \lcite[Thm.\ 38.2]{Re75}.
\end{proof}
\end{proposition}
\begin{remark}
Cancellation is not always true, see \cite{shioda77::some+abelian}. Its
characteristic zero example involves a non-hereditary order. However, the purely
algebraic consequences of one of its positive characteristic example implies that
there is an example where $R$ is a maximal order in a definite quaternion
algebra over $\Q$.
\end{remark}
Recall that the Kummer sequence (resp.\ the exponential sequence) gives a group
homomorphism from the Brauer group of a variety $X$ (resp.\ over $\C$) to
$H^3(X,\widehat \Z)$ (resp.\ $H^3(X,\Z)$). The image of an element of the Brauer
group under this map will be called its \Definition{characteristic class}.
%;;%%Proposition BS non-triviality
\begin{proposition}\label{BS non-triviality}
Suppose $\k$ is algebraically closed. For every integer $n\ge 1$ for which $n+1$
is not a power of the characteristic of $\k$ there is a smooth and projective
variety $X$ and a projective space fibration $P \to X$ of relative dimension $n$
and with non-trivial characteristic class.
\begin{proof}
Suppose that we have a finite group $G$ of order prime to the characteristic if
the latter is positive acting freely on a variety $Y$ and a projective
representation \map{\pi}{G}{\PGL_{n+1}(\k)}. We then get a projective space
bundle $P:=Y\times_G\P^n \to Y/G=:X$ and it is clear that its characteristic
class is the inverse image under the classifying map $X \to \Bclass G$ of $Y \to
X$ of the characteristic class of $\pi$. The latter is defined as the image of
$\pi$ under the composite $H^1(G,\PGL_{n+1}(\k)) \to H^2(G,\k^*) \to
H^2(G,\bk^*)=H^2(G,\Q/\Z)=H^3(G,\Z)$.

Now I claim that for every positive integer $K$ there is a smooth and projective
$Y$ as above such that the induced map $H^i(G,\Z) \to H^i(X,\widehat\Z)$ is
injective for $i \le K$. Indeed, the Godeaux-type construction of Serre (cf.,
\cite[\S20]{serre58::sur}) starts with a linear representation $V$ of $G$ whose
associated projective representation is faithful and constructs a smooth complete
intersection $Y$ in the projective space $\P(V)$ of $V$ on which $G$ acts
freely. The maximal possible dimension of $Y$ depends on the codimension of the
points of $\P(V)$ with non-trivial stabiliser. That codimension may be made as
large as possible by choosing a suitable $V$ (for instance by replacing a given
$V$ by a high power of it) and therefore we may find $Y$'s of arbitrarily high
dimension. Consider first the stack quotient $[\P(V)/G]$. The natural map
$[\P(V)/G] \to \Bclass G$ makes $[\P(V)/G]$ a projective bundle over $\Bclass G$
associated to the vector bundle $[V/G] \to \Bclass G$. This implies that if $\xi
\in H^2([\P(V)/G],\widehat\Z)$ is the first Chern class of the $G$-linearised
line bundle $\sO(1)$, then $H^*([\P(V)/G],\widehat)$ is the free
$H^*(G,\widehat\Z)$-module on the $\xi^i$, $i=0,\dots,n-1$ (cf.,
\cite[1.9]{grothendieck68::class+chern}). (There may be some doubt about the
$\Z_p$-part of this when $p=\Char\k$. It turns out to be true but is not used in
the subsequent arguments.) We now have a morphism $X=Y/G=[Y/G] \to [\P^n/G]$
induced by the inclusion $Y \subseteq \P^n$ giving rise to a comparison map from
the spectral sequence $H^i(G,H^j(\P^n,\widehat\Z)) \implies
H^{i+j}([\P^n/G],\widehat\Z)$ to the corresponding one for $[Y/G]$. By the
Lefschetz theorem the map $H^i(\P^n,\widehat\Z) \to H^i(Y,\widehat\Z)$ is an
isomorphism for $i<\dim Y$ which by the comparison map of spectral sequences
implies that $H^i([\P^n/G],\widehat\Z) \to H^i(X,\widehat\Z)$ is an isomorphism
for $i < \dim Y$. In particular, the composite $H^i(G,\widehat\Z) \to
H^i(X,\widehat\Z)$ is injective in the same range. As we can make $\dim Y$
arbitrarily large we have proven our claim.

Let now $1 \to \k^* \to H \to G \to 1$ be a non-trivial central extension of $G$
and assume that we have a linear representation $H \to \GL_{n+1}(\k)$ of $H$
which takes $\k^*$ to scalar matrices. This gives us an induced projective
representation $G \to \PGL_{n+1}$. I claim that the characteristic class of this
representation is equal to the class of the central extension $H$. Indeed, while
the most conceptual would no doubt be ``gerbic'', the simplest is probably the
easy computation by cocycles. In particular it is non-trivial as the central
extension is. Now, choose a $Y$ as above with $K\ge 3$. The pullback of $\alpha$
to $X$ represents the projective bundle $Y\times_G\P^n \to Y/G$ and by
assumption it is non-zero. If $n+1$ is not divisible by the characteristic of
$\k$ we can let $H$ be the $g=1$ level $n+1$ theta group so that
$G=\Z/(n+1)\Z\times\Z/(n+1)\Z$ and the theta group has an $n+1$-dimensional
representation of the required kind. If $n+1$ is not a power of the
characteristic we may let $d$ be any non-trivial divisor of $n+1$ not divisible
by the characteristic and let $H$ be the level $d$ theta group and let the
$n+1$-dimensional representation be $(n+1)/d$ copies of the standard
representation.
\end{proof}
\end{proposition}
\begin{remark}
\part As we only need $K\ge 3$ (with the notations of the proof) the proof
actually gives us $3$-dimensional examples.

\part There are $2$-dimensional examples at least for some values of $n$. One
such example is an Enriques surface. In that case Poincar\'e duality shows that
the torsion subgroup of $H^3(X,\widehat \Z)$ is equal to $\Z/2$ and it comes
from the Brauer group (as the cohomological Brauer group is equal to the Brauer
group). I don't know however which $n$ can be chosen in this case (presumably
$n=1$ is possible).
\end{remark}
%;;%%Lemma Fibration formula
\begin{lemma}\label{Fibration formula}
Let $G$ be a semi-simple group over an algebraically closed field $\k$.

\part $\{G\}$ is invertible in $\Kscheme{\k}'$.

\part Assume that $F$ is a $G$-space and $P \to X$ is a $G$-torsor (of
spaces). If we for the two $G$-torsors $P \to X$ and $P\times F \to P\times_GF$
have the equalities $\{P\}=\{G\}\{X\}$ and $\{P\times
F\}=\{G\}\{P\times_GF\}$, then we have $\{P\times_GF\}=\{F\}\{X\}$.
\begin{proof}
If we let $B$ be a Borel subgroup of $G$ we have a $B$-torsor $G \to G/B$ and as
$B$ is a successive extension of $\mul$'s and $\add$'s and thus $B \in \Zar$ so
that $\{G\}=\{B\}\{G/B\}$. Now, again by the fact $B$ is a successive extension
$\mul$'s and $\add$'s we have that $\{B\}$ is an integer polynomial in $\Lclass$
and as $G/B$ has a cell decomposition so does $\{G/B\}$. This polynomial is
determined by the fact that $\chi_c(\{G\})$ is the same polynomial in
$q:=\chi_c(\Lclass)$ and it is well-known that this polynomial is a power of $q$
times a product of $q^n-1$'s which shows that $\{G\}$ is indeed invertible in
$\Kscheme{\k}'$.

As for the second part we have, using the assumed relations,
\begin{displaymath}
\{G\}\{F\}\{X\}=\{F\}\{P\}=\{P\times F\}=\{G\}\{P\times_GF\}
\end{displaymath}
and by the first part we can cancel $\{G\}$.
\end{proof}
\end{lemma}
%;;%%Proposition PSL2 example
\begin{proposition}
Assume that $\k$ is algebraically closed. There are $\PSL_{n}$-torsors, for
any $n\ge 1$ not a power of the characteristic of $\k$, $T
\to X$ over smooth and projective varieties $X$ such that $\{T\}\ne
\{\PSL_n\}\{X\} \in \Kschemecompl{\k}$. In particular the map $\Kschemecompl{\k}
\to \Kschemecompl[\PSL_n]{\k}$ is not injective and similarly for its variants.
\begin{proof}
By Proposition \ref{BS non-triviality} there is a $\P^{n-1}$-fibration
\map{\pi}{P}{X}, $X$ smooth and proper, whose characteristic class is
non-trivial. The pullback of $P \to X$ along $P \to X$ has a section and hence
has a reduction of the structure group to $\GL_{n+1}$ and in particular the
class $\alpha$ pulls back to $0$ in $H^3(P,\Z_p)$. Now, we have that the
$R^{i}\pi_*\Z_p$ are constant as they become constant after pullback along $\pi$
which has connected fibres.  As also $R^{2i+1}\pi_*\Z_p=0$, the Leray spectral
sequence for $\pi$ gives an exact sequence
\begin{displaymath}
H^3(X,\Z_p) \to H^3(P,\Z_p) \to H^1(X,R^2\pi_*\Z_p)=H^1(X,\Z_p).
\end{displaymath}
Now, tensored with $\Q$, the spectral sequence degenerates (for instance for
reasons of weight) so that the kernel of $H^3(X,\Z_p) \to H^3(P,\Z_p)$ is a
torsion group and furthermore $H^1(X,\Z_p)$ is torsion free. Together this shows
that the torsion subgroup of $H^3(X,\Z_p)$ maps surjectively onto that of
$H^3(P,\Z_p)$. However, we have seen that $\alpha$ goes to zero so the
conclusion is that the torsion subgroup of $H^3(P,\Z_p)$ is strictly smaller
than that of $H^3(X,\Z_p)$. This implies that we do \emph{not} have that
$\{P\}=\{\P^n\}\{X\} \in \Kschemecompl{\k}$ as otherwise
$H^\bullet(\chi_G(P))=(1-s^{n+1})/(1-s)H^\bullet(X)$, where $s := t^{2}$, and in
particular $\{H^3(P,\Z_p)\}$ would be equal to
$\{H^3(X,\Z_p)\}+\{H^1(X,\Z_p)\}$. This would give that the torsion subgroup of
$H^3(P,\Z_p)$ would be equal to that of $H^3(X,\Z_p)$ which we have seen is not
the case. 

Following Lemma \ref{Fibration formula} one of the two $\PSL_{n+1}$ torsors, $T
\to X$ associated to $P \to X$ and $T\times \P^n \to T\times_{\PSL_{n+1}}
\P^n=P$ do not fulfil the multiplicativity relation.
\end{proof}
\end{proposition}
\begin{remark}
\part Lemma \ref{Fibration formula} is necessary only to pinpoint the torsors for
which multiplicativity fails. If one doesn't care about that one could, by way
of contradiction, assume multiplicativity for all $\PSL_{n+1}$-torsors which
implies that $\Kschemecompl{\k} \to \Kscheme[\PSL_{n+1}]{\k}$ is an
isomorphism. We can then apply (\ref{Fibration}).

\part Similar results are true for every semi-simple special algebraic group and
will be shown elsewhere.
\end{remark}
\end{section}
%;;%Section A sample computation
\begin{section}{A sample computation}

We shall now give an example in the spirit of \cite{bergstroem06::point} (though
the cases covered by Bergstr\"om are more complicated than the one we shall
discuss). The general technique is the following: Present the stack whose class
we want to compute as the stack quotient of a linear group acting (linearly) on
an open subset of a linear space (or possibly projective space, though that
would typically give us $\PGL_n$ as linear group which does have torsors not
trivial in the Zariski topology). The class of the stack quotient for the action
of the group on the full linear space is equal to a power of $\Lclass$ divided
by the class of the group. Hence, we may look instead at the complement of the
open subset. That complement can then be stratified (typically by imposing
various kinds of singular sets) and for each stratum the data specifying it can
be put in standard position (which just means a reduction of the structure
group). A stratum with data in fixed position can then be embedded as an open
subset of a linear space with group action and one is again reduced to
considering the complement. Unfortunately the number of strata that one must
consider very quickly becomes large so I shall be somewhat sketchy as my aim is
only to illustrate the technique.
\begin{example}
We shall need to compute $\{\Bclass\Sigma_2\}$ and $\{\Bclass\Sigma_3\}$. In
\cite{ekedahl08::stack} we shall prove that in general $\{\Bclass\Sigma_n\}=1$
for all $n$ but these cases can be done by hand. As we are going to make that
restriction anyway we assume that the characteristic is different
from $2$ and $3$. For $\Sigma_2\iso\{\pm1\}$ we can embed it in $\mul\in \Zar$
and hence by Proposition \ref{Base formulas} we get $\{\Bclass
\Sigma_2\}=\{\mul/\pm1\}\{\Bclass\mul\}=\{\mul\}\{\Bclass\mul\}$ and
$1=\{\mul\}\{\Bclass\mul\}$ which gives $\{\Bclass\Sigma_2\}=1$. For
$\Bclass\Sigma_3$ we have that the standard representation divided by the
trivial subrepresentation gives an embedding $\Sigma_3 \subset \GL_2\in\Zar$ and
hence by Proposition \ref{Base formulas} we get that
$\{[\P^1/\Sigma_3]\}=\{\P^1\}\{\Bclass\Sigma_3\}$. Now, there are three points
of $\P^1$ with non-trivial $\Sigma_3$-fixed points and if $U$ is the complement
of them we have $[U/\Sigma_3]=U/\Sigma_3$ and $U/\Sigma_3$ is $\P^1$ minus one
point. Hence if $X:=\P^1\setminus U$ we have
$\{[\P^1/\Sigma_3]\}=\{U/\Sigma_3\}+\{[X/\Sigma_3]\}=\{\P^1\}-1+\{[X/\Sigma_3]\}$.
As $X$ consists of three points permuted by $\Sigma_3$ according to the standard
permutation representation and hence $[X/\Sigma_3]\iso \Bclass\Sigma_2$ and by
what we just proved $\{[X/\Sigma_3]\}=1$ and thus $\{[\P^1/\Sigma_3]\}=\{\P^1\}$
and thus we get $\{\P^1\}=\{\P^1\}\{\Bclass\Sigma_3\}$ and as
$\{\P^1\}=\Lclass-1$ it is invertible so we get $\{\Bclass\Sigma_3\}=1$.
\end{example}
Let $\sM$ be the algebraic stack of smooth and proper genus $1$ curves with a
polarisation of degree $3$. Here \emph{polarised} means a curve together with a
line bundle. To the universal pair $(\sC \to \sM,\sL)$ we associate the rank $3$
vector bundle $\pi_*\sL$ (where $\pi$ is the projection $\sC \to \sM$). The
frame bundle $E$ of $\pi_*\sL$ is spatial, in fact the space of non-singular
cubic forms in three variables, and we have that $\sM$ is the stack quotient
$[E/\GL_3]$. Now, let $V$ be the space of cubic forms and $V'$ the closed
subspace of singular cubic forms. Thus we have
\begin{displaymath}
\Lclass^{10}\{\GL_3\}^{-1}= \{[V/\GL_3]\} = \{[V'/\GL_3]\} + \{\sM\}
\end{displaymath}
and as the left hand side is known, we focus attention on $\{[V'/\GL_3]\}$. To
continue we assume, for simplicity, that $\Char\k\ne2,3$. Inside of $V'$ we have
the open subset $V'_1$ of cubics with exactly one singular point. In turn $V'_1$
can be divided up into to the nodal part $V'_n$, the open subset where the
singular locus is reduced and the cuspidal part $V'_c$, where also the second
derivatives of the cubic vanish at the singular point. In both cases one may fix
a point, $(0\hco 1\hco 0)$, and let $U'_n$ and $U'_c$ be the set of cubics
having the fixed point as unique singular point which is nodal resp.\
cuspidal. In the nodal case the singular locus gives a section and from that it
easily follows that the natural map $U'_n \to [V'_n/\GL_3]$ is smooth and
surjective and from that one gets an isomorphism $[U'_n/G] \riso [V'_n/\GL_3]$,
where $G$ is the stabiliser in $\GL_3$ of the fixed point. Similarly, in the
cuspidal case, the vanishing of the cubic, the first and the second derivatives
defines a section giving an isomorphism $[U'_c/G] \riso [V'_c/\GL_3]$. These
isomorphisms of stacks of course give equalities $\{[U'_n/G]\}=\{[V'_n/\GL_3]\}$
and $\{[U'_c/G]\}=\{[V'_c/\GL_3]\}$ and we can then make a reassembly: We let
$U'_1$ be the space of cubics with a unique singularity at the fixed point and
we then get
\begin{displaymath}
\{[U'_1/G]\} = \{[U'_n/G]\}+\{[U'_c/G]\} = \{[V'_n/\GL_3]\}+\{[V'_c/\GL_3]\}=\{[V'_1/\GL_3]\}.
\end{displaymath}
Continuing we let $U$ be the linear space of cubics with a singular point at
$(0\hco 1\hco 0)$ so that $U'_1$ is an open subset of $U$. Again we get that
$\{[U/G]\}=\Lclass^{7}\{G\}^{-1}$, as $G\in \Zar$ being the extension of groups
that are, (and as scheme $G\iso
\GL_2\times\GL_1\times\A^2$ so that $\{G\}=\{\GL_2\}\{\GL_1\}\Lclass^2$). We are
therefore reduced to considering the complement $U''$ of $U'_1$ in $U$. Going
back to an earlier reduction we also need to consider the complement $V''$ of
$V'_1$ in $V'$. This leads to the following cases:
\begin{itemize}
\item We have a stratum $W$ in $V$ consisting (geometrically) of three distinct
lines not all passing through a point. As $\{[W/\GL_3]\}=\{W\}\{\GL_3\}^{-1}$ it
is enough to compute $\{W\}$. Letting $\overline W$ be the quotient of $W$ by
the action by the scalars, the quotient map $W \to \overline W$ is a
$\mul$-torsor so we have $\{W\}=(\Lclass-1)\{\overline W\}$ and it will be
enough to compute $\{\overline W\}$.  We have a map from $W'$, the space of
linearly independent triples of vectors (a.k.a.\ invertible matrices) in
$(\A^3)^3$, to $W$, given by thinking of the vectors as linear forms and
multiplying them together. We can act by the semi-direct product
$H:=\mul^3\ltimes\Sigma_3$ on $W'$, multiplying each vector by an invertible
scalar and permuting the three vectors. The action is free and we have
$\overline{W}=W'/H$. The action is also isomorphic to left multiplication of $H$
on $\GL_3$, where $H$ is identified with the normaliser of the group of diagonal
matrices in $\GL_3$. Hence we get from (\ref{Subgroup formula}) that
$\{W\}=\{\Bclass H\}\{\GL_3\}$ and it remains to compute $\{\Bclass H\}$. For
that we consider the natural linear action of $H$ on $V=\k^3$. As usual we have
that $\{[V/H]\}=\Lclass^3\{\Bclass H\}$. The diagonal matrices of $H$ act freely
and transitively on the open subset $V_0$ of $V$ defined by all coordinates
being non-zero. This gives that $[V_0/H]\iso\Bclass\Sigma_3$ and we have already
noticed that $\{\Bclass\Sigma_3\}=1$. We then look at the stratum $V_1$ where
exactly one coordinate is zero. Again $H$ acts transitively so that $[V_1/H]$ is
isomorphic to the classifying stack of the stabiliser of $(1,1,0)$ say. This
stabiliser is isomorphic to $\Sigma_2\times\mul$ so that
$\{[V_1/H]\}=\{\Bclass\Sigma_2\}\{\Bclass\mul\}$. We have however shown that
$\{\Bclass\Sigma_2\}=1$. Similarly, for $V_2$ the locus where exactly two
coordinates are zero we get $[V_2/H]\iso\Bclass\mul^2$ and we finally get a
contribution $\{\Bclass H\}$ from $[\{0\}/H]$. Moving over this last term to the
left hand side of $\{[V/H]\}=\Lclass^3\{\Bclass H\}$ we get a formula for
$(\Lclass^3-1)\{\Bclass H\}$ and as $\Lclass^3-1$ is invertible in $\Kstack{\k}$
we get a formula for $\{\Bclass H\}$, more precisely we get that $\{\Bclass
H\}=(\Lclass-1)^3$.

%% It is clear that the action is free and
%% that $W=W'/H$. The action of $H$ on $W'$ is the restriction to $W'$ of a linear
%% action on $(\A^3)^3$ and we have $\{[(\A^3)^3/H]\}=\Lclass^9\{\Bclass H\}$ and,
%% provided we can compute $\{\Bclass H\}$ we are reduced to considering the
%% complement of $W'$ in $(\A^3)^3$. The complement is divided into strata where
%% the span of the vectors have fixed dimension. The rank $2$-locus is given by a
%% point $\ell$ and a line $U$ in $\P^2$ together with an isomorphism $\k^3/\ell
%% \riso U$. Hence we are in the situation of (\ref{Isotorsor}) and we get that the
%% class of that locus equals $\{\GL_2\}\{\P^2\}^2$. Similarly, the rank $1$-locus
%% has class $\{\GL_1\}\{\P^2\}^2$ and the rank $0$-locus class $\}\{\P^2\}^2$.

\item For $U$ we need to consider three distinct lines two of which passes
through a fixed point. The analysis is altogether similar to the previous one.

\item We have one stratum consisting of a line and a smooth quadric (resp.\ with
a singular point at a fixed point). We can place the line in standard position
and compare the space of non-singular quadrics (in the second case passing
through a fixed point) with the linear space of all quadrics. The space of
singular quadrics is then easily analysed.

\item We have the case of a double line and another distinct line and a triple
line. The analysis is somewhat simplified by the fact that we work in
characteristic different from $2$ and $3$ which means that the map $\ell \mapsto
\ell^2$ from lines to quadrics is an immersion (and similarly for cubes).
\end{itemize}
One may indeed go through all the calculations to get an explicit formula for
$\{\sM\}$. One could also argue as follows: We can check quite easily that each
contribution to the result is a rational function in $\Lclass$ in
$\Kstack{\k}=\Kscheme{\k}'$. We can the use the fact that $\chi_c$ is injective
on such rational functions, indeed it is clear that if $\Sigma_{i\le N}a_iq^i=0$
where the $a_i$ are integers and $q=\chi_c(\Lclass)$ then $a_i=0$, and hence it
is enough to determine $\chi_c(\{\sM\})\in \Kcpol{\k}$ as a rational function in
$\chi_c(\Lclass)$. By Theorem \ref{Real Euler characteristic} this can be done
by taking the ordinary Euler characteristic with compact support. Furthermore,
$\sM$ is a smooth stack so the cohomology of compact support is dual to the
ordinary cohomology. Now, if $\sM_{0,1}$ is the stack of smooth proper genus $1$
curves with a distinguished point we have a map $\sM_{0,1} \to \sM$ given by
$(E,p) \mapsto (E,\sO(3p))$. It is clear that it induces an isomorphism on
cohomology with characteristic zero coefficients. Hence, we get that
$\chi_c(\{\sM\})=\chi_c(\{\sM_{0,1}\})=\chi_c(\Lclass)$ and we conclude that
$\{\sM\}=\Lclass$.
\end{section}
\bibliography{preamble,abbrevs,alggeom,algebra,ekedahl}
\bibliographystyle{pretex}
\end{document}